\documentclass[11pt]{amsart}
\usepackage[utf8]{inputenc}
\usepackage{changepage}
\usepackage{subfiles}

\usepackage{amsmath}
\usepackage{amssymb}
\usepackage{amsthm}
\usepackage{thmtools}
\usepackage{thm-restate}
\usepackage[margin=1in]{geometry}
\usepackage{cleveref}
\usepackage{ytableau}

\usepackage{biblatex}
\renewbibmacro{in:}{}

\numberwithin{equation}{section}

\usepackage{cleveref}
\usepackage{float}
\usepackage{xcolor}
\usepackage[utf8]{inputenc} 
\usepackage[T1]{fontenc}
\usepackage{comment}
\usepackage{lscape}

\usepackage{enumitem}
\usepackage{amsmath}
\usepackage{tikz}
\usepackage{multicol}
\usepackage[center]{caption}

\usetikzlibrary{decorations.pathreplacing}
 \newtheorem{thm}{Theorem}
 
 \newtheorem{lem}[thm]{Lemma}
 
 \theoremstyle{definition}
 
 \theoremstyle{remark}

 \numberwithin{equation}{section}
\newcommand{\Cross}{$\mathbin{\tikz [x=1.4ex,y=1.4ex,line width=.2ex] \draw (0,0) -- (1,1) (0,1) -- (1,0);}$}%
\usetikzlibrary{shapes.arrows}

\title{Counting crucial permutations with respect to monotone patterns}
\author{Yunseo Choi}

\begin{document}

\maketitle

\begin{abstract}
Recently, Avgustinovich, Kitaev, and Taranenko defined five types of $(k, \ell)-$crucial permutations, which are maximal permutations that do not contain an increasing subsequence of length $k$ or a decreasing subsequence of length $\ell$. Further, Avgustinovich, Kitaev, and Taranenko began the enumeration of the $(k, \ell)-$crucial permutations of the minimal length and the next minimal length and the $(k, 3)-$crucial permutations of all lengths for each of the five types of $(k,\ell)-$crucial permutations. In this paper, we complete the enumeration that Avgustinovich, Kitaev, and Taranenko began. 
\end{abstract}

\section{Introduction} 
\label{intro}

The classical Erd\H{o}s-Szekeres theorem \cite{erdos_szekeres} from 1935 states that a permutation of length $(k-1)(\ell-1)+1$ must contain an increasing subsequence of length $k$ or a decreasing subsequence of length $\ell$. Accordingly, the permutations of length $(k-1)(\ell-1)$ that do not contain an increasing subsequence of length $k$ or a decreasing subsequence of length $\ell$ are called the \textit{$(k, \ell)$-Erd\H{o}s-Szekeres extremal permutations} \cite{vatter}. The number of $(k, \ell)$-Erd\H{o}s-Szekeres extremal permutations for $k \geq \ell$, which follows from the Robinson-Schensted-Knuth (RSK) correspondence  \cite{knuth, robinson, schensted} and the hook-length formula \cite{frt}, is given by 
\begin{equation} \label{Erd\H{o}s-szekeres}
    \left(\frac{((k-1)(\ell-1))!}{1^1 2^2 \cdots (\ell-1)^{\ell-1} \ell^{\ell-1} \cdots (k-1)^{\ell-1} k^{\ell-2} \cdots (k+\ell-3)^{1}}\right)^2. 
\end{equation} The formula was first presented for $k = \ell$ by Stanley \cite{stanley} in 1969 and its generality by Avgustinovich, Kitaev, and Taranenko \cite{akt} in 2022. 

In a larger context, \textit{crucial objects} \cite{crucial_2, crucial_1} in combinatorics are objects that satisfy a set of conditions but cease to satisfy the conditions when extended in any allowed way. Recently, Avgustinovich, Kitaev, and Taranenko \cite{akt} defined \textit{$(k, \ell)$-right-crucial permutations} to be the permutations that are crucial when the set of conditions is avoid increasing subsequences of length $k$ and decreasing subsequences of length $\ell$ and the allowed extensions are to the right. More formally, let $\pi= \pi_1 \pi_2 \cdots \pi_{n} \in S_{n}$. The \textit{$i^\text{th}$ extension of $\pi$ to the right (respectively, left)} is the permutation $f(\pi_1) f(\pi_2) \cdots f(\pi_{n}) i$ (respectively, $i f(\pi_1) f(\pi_2) \cdots f(\pi_{n})$), where $$ f(j) = \begin{cases} 
      j & \text{if } j < i \\
      j+1 & \text{if } j \geq i.
   \end{cases}$$

In addition to the $(k, \ell)$-right-crucial permutations, Avgustinovich, Kitaev, and Taranenko \cite{akt} defined four other types of \textit{$(k, \ell)$-}crucial permutations. Each type of $(k, \ell)$-crucial permutations shares the set of conditions with the $(k, \ell)$-right-crucial permutations but differ in the extensions that are allowed. The following list showcases the five types of $(k, \ell)$-crucial permutations that Avgustinovich, Kitaev, and Taranenko defined alongside their allowed extensions. The \textit{$i^\text{th}$ extension of $\pi$ from above (respectively, below)} is the permutation $\pi_{1} \pi_{2} \cdots \pi_{i-1} (n+1) \pi_i \cdots \pi_n$ (respectively, $(\pi_{1}+1) (\pi_{2}+1) \cdots (\pi_{i-1}+1) 1 (\pi_i+1) \cdots (\pi_n+1)$). As noted by Avgustinovich, Kitaev, and Taranenko (see Lemma 11 of \cite{akt}), the list below is the exhaustive list of $(k, \ell)$-crucial permutations up to bijections when the set of extensions is some subset of the extensions to the left, to the right, from above, or from below. 
\begin{enumerate}
    \item $(k, \ell)$-right-crucial: any extension to the right 
    \item $(k, \ell)$-bicrucial: any extension to the right or the left
    \item $(k, \ell)$-top-right-crucial: any extension to the right or from above
    \item $(k,\ell)$-tricrucial: any extension to the right or the left or from above
    \item $(k, \ell)$-quadrocrucial: any extension to the right or the left or from above or below
\end{enumerate}

Now, a rephrasing of the Erd\H{o}s-Szekeres theorem is that $(k, \ell)$-Erd\H{o}s-Szekeres extremal permutations are $(k, \ell)$-crucial in any of the five ways. But whereas the $(k, \ell)$-Erd\H{o}s-Szekeres extremal permutations must have length $(k-1)(\ell-1)$, Avgustinovich, Kitaev, and Taranenko (see Lemma 4 and Theorem 10 of \cite{akt}) proved that the length of a $(k, \ell)$-(right, top-right)-crucial can be as short as $k+\ell-3$. Similarly, they proved that a $(k, \ell)$-(bi, tri, quadro)crucial permutation can be as short as length $k+\ell+ \mathrm{min}\{k, \ell\}-5$ (see Theorem 9, Theorem 32, and Theorem 36 of \cite{akt}). Following Avgustinovich, Kitaev, and Taranenko, we refer to the $(k, \ell)$-crucial permutations of the shortest (respectively, next shortest) length as the \textit{minimal (respectively, next minimal) $(k, \ell)$-crucial permutations}.

In the spite of \Cref{Erd\H{o}s-szekeres}, which enumerates the $(k, \ell)$-crucial permutations of the longest length, Avgustinovich, Kitaev, and Taranenko \cite{akt} began the enumeration of special classes of $(k, \ell)$-crucial permutations: namely, the minimal and next minimal $(k, \ell)$-crucial permutations and the $(k, 3)$-crucial permutations of all lengths. In this paper, we complete the enumeration that they began for all five types of $(k, \ell)$-crucial permutations. As do Avgustinovich, Kitaev, and Taranenko, we base our enumeration on the RSK correspondence and the hook-length formula. 

The rest of the paper is organized as follows. In \Cref{prelim}, we establish the preliminaries. In \Cref{right}, we count the next minimal $(k, \ell)$-right-crucial permutations. In \Cref{bi}, we count the next minimal $(k, \ell)$-bicrucial permutations and the $(k, 3)$-bicrucial permutations. In \Cref{tri}, we count the minimal $(k, \ell)$-minimal tricrucial permutations, the next minimal $(k, \ell)$-tricrucial permutations, and the $(k, 3)$-tricrucial permutations. Finally, in \Cref{quadro}, we count the next minimal $(k, \ell)$-quadrocrucial permutations and the $(k, 3)$-quadrocrucial permutations. \Cref{summary} summarizes our findings.

\begin{footnotesize}
\begin{center}
\begin{table}[htpb]
    \centering
    \begin{tabular}{c|c|c|c}
         & \begin{tabular}{@{}c@{}} enumeration of minimal \\ $(k,\ell)$-crucial permutations\end{tabular} & \begin{tabular}{@{}c@{}} enumeration of the next minimal \\ $(k,\ell)$-crucial permutations\end{tabular} & \begin{tabular}{@{}c@{}} enumeration of \\ $(k,3)$-crucial permutations\end{tabular}  \\ \hline \hline
         right-crucial & Corollary 8 of \cite{akt} & \Cref{right_crucial_next} & Corollary 13 of \cite{akt} \\
         top-right-crucial & Theorem 10 of \cite{akt} & Theorem 29 of \cite{akt} & Corollary 27 of \cite{akt} \\
         bicrucial & Corollary 21 of \cite{akt} & \Cref{bicrucial_next} & \Cref{bicrucial_3} \\
         tricrucial & \Cref{tricrucial_min} & \Cref{next_min_tri} & \Cref{tricrucial_3} \\
         quadrocrucial & Theorem 37 of \cite{akt} & \Cref{next_min_quadro} & \Cref{quadrocrucial_3} 
    \end{tabular}
    \caption{A summary of our results}
     \label{summary}
\end{table}
\end{center}
\end{footnotesize}

\section{Preliminaries}
\label{prelim}
A \textit{partition} $\lambda = (\lambda_1, \ldots, \lambda_k)$ of $n$ is a weakly decreasing sequence of non-negative integers such that $\lambda_1 + \cdots + \lambda_k = n$. A \textit{Young diagram of shape $\lambda$} is a left justified diagram of boxes with $\lambda_i$ boxes in the $i^{\text{th}}$ row for $i \in \{1, 2, \ldots, k\}$. A \textit{skew Young diagram of shape $\lambda$} is a diagram of boxes with $\lambda_i$ boxes in the $i^{\text{th}}$ row that is not necessarily left and top justified. We say that $\lambda_{1} + \lambda_{2} + \cdots + \lambda_{k}= n$ is the \textit{size} of a (skew) Young diagram. A (skew) Young diagram is \textit{filled} if each box contains exactly one number, and the numbers that appear in the boxes are distinct. A \textit{(skew) standard Young tableau of shape $\lambda$} is a (skew) Young diagram of shape $\lambda$ filled with numbers $1, 2, \ldots, n$ such that every row is increasing from the left to right and every column is increasing from the top to bottom. For a (skew) standard Young tableau $T$, let $t_{i, j}$ (respectively, $b_{i, j}$) be the entry (respectively, the box) at the $i^{\text{th}}$ row and the $j^{\text{th}}$ column of $T$. We define the \textit{complement} $T^{c}$ of a (skew) standard Young tableau $T$ to be the (skew) Young diagram that shares the shape with $T$ and satisfies $t^{c}_{i, j} = n+1-t_{i,j}$ for all entries $t^{c}_{i,j}$ in $T^{c}$. We define the \textit{transpose} $T^{t}$ of a (skew) standard Young tableau $T$ to be the (skew) Young diagram that shares the shape with $T$ after a reflection across its main diagonal and satisfies $t^{t}_{i, j} = t_{j,i}$ for all entries $t^{t}_{i,j}$ in $T^{t}$.  

The Robinson-Schensted-Knuth (RSK) correspondence \cite{knuth, robinson, schensted} is a bijection between the permutations of length $n$ and the pairs of standard Young tableaux $(P, Q)$ that have the same shape $\lambda$. The bijection from $S_n$ to $(P, Q)$ sends $\pi$ to $(P_{n}(\pi), Q_{n}(\pi))$. The construction of $P_{n}(\pi)$ builds off of \textit{row insertion $T \leftarrow i$} for a tableau $T$ and an element $i$. Let $T_1$ be the first row of $T$ and $T_2$ be the filled Young diagram of $T$ with $T_1$ removed. The procedure $T \leftarrow i$ is defined recursively as follows. 
\begin{enumerate}
    \item If $T$ is $\emptyset$ or the largest entry in $T_1$ is less than $i$, then insert $i$ at the right end of $T_1$.
    \item Otherwise, let $j$ be the smallest element in $T_1$ that is greater than $i$. Replace $j$ with $i$ in $T_1$, and let the remaining rows of be given by $T_2 \leftarrow j$. 
\end{enumerate}

$P_{i}(\pi)$ is now defined as follows. For $\pi = \pi_1 \pi_2 \cdots \pi_n \in S_n$, we let $$P_{i}(\pi) = (\cdots((\emptyset \leftarrow \pi_1) \leftarrow \pi_2) \cdots \leftarrow \pi_{i-1}) \leftarrow \pi_i.$$

Next, we define $Q_{n}(\pi)$. Let $Q_0 = \emptyset$, and recursively let $Q_i(\pi)$ be the tableau obtained by inserting $i$ to $Q_{i-1}(\pi)$ such that $Q_i(\pi)$ and $P_{i}(\pi)$ have the same shape and the position of the entries of $Q_{i-1}(\pi)$ remain invariant from $Q_{i-1}(\pi)$ to $Q_{i}(\pi)$. 

For more details on the RSK correspondence, see, for example, \cite{rsk}. As we proceed, let $p_{i, j}$ (respectively, $q_{i, j}$) be the entry at the $i^{\text{th}}$ row and the $j^{\text{th}}$ column of $P_{n}(\pi)$ (respectively, $Q_{n}(\pi)$). Furthermore, as usual, we let $\delta_{a,b}$ be the Kronecker delta, where
$$
\delta_{a,b}:=\begin{cases}
			1, & \text{if } a = b \text{ and}\\
            0 & \text{otherwise}.
		 \end{cases}
$$ In addition, we let $\binom{a}{b_{1},b_{2}, \cdots, b_{m}}$ denote the multinomial coefficient 
$$
 \binom{a}{b_{1}, b_{2}, \ldots, b_{m}} := \frac{a!}{b_{1}! b_{2}! \ldots b_{m}!}. 
$$

\section{Right-Crucial Permutations}
\label{right}
In this section, we enumerate the next minimal $(k, \ell)$-right-crucial permutations. We begin by citing a result by Avgustinovich, Kitaev, and Taranenko \cite{akt}, which characterizes the RSK correspondence $(P_{n}(\pi), Q_{n}(\pi))$ of $(k, \ell)$-right-crucial permutations $\pi$.

\begin{thm}[Avgustinovich, Kitaev, and Taranenko \cite{akt}] \label{right_crucial_char}
The permutation $\pi \in S_{n}$ is $(k,\ell)$-right-crucial if and only if:
\begin{enumerate}
    \item the number of columns in $P_{n}(\pi)$ (and $Q_{n}(\pi)$) is $k-1$;
    \item the number of rows in $P_{n}(\pi)$ (and $Q_{n}(\pi)$) is $\ell-1$;
    \item $P_{n}(\pi)$ contains an increasing sequence of elements $p_{1, k-1}, p_{2, j_2}, \ldots, p_{\ell-1, j_{\ell-1}}$ for some $j_2,$ $ \ldots,$ $j_{\ell-1}$. 
\end{enumerate}
\end{thm}

\subsection{Next minimal right-crucial permutations}

We cite another result, which is implied from Theorem 29 of Avgustinovich, Kitaev, and Taranenko \cite{akt}. 

\begin{lem}[Avgustinovich, Kitaev, and Taranenko \cite{akt}] \label{lem_for_next_2}
The number of standard Young tableaux $T$ of shape $(k-1, 2, \underbrace{1, \ldots, 1}_{\ell-3})$ that satisfy $t_{1, k-1} < t_{2, j} < t_{3, 1}$ for some $j$ is $k + \ell-4$. 
\end{lem}

Our enumeration of the next minimal $(k, \ell)$-right-crucial permutations now follows from the hook-length formula, \Cref{right_crucial_char}, and \Cref{lem_for_next_2}. 

\begin{thm} \label{right_crucial_next}
For $k, \ell \geq 3$, the number of next minimal $(k, \ell)$-right-crucial permutations is 
$$
\frac{(k-2)(\ell-2)(k+\ell-4)}{k+\ell-3} {l+k-2 \choose l-1}. 
$$
\end{thm}
\begin{proof}
By Theorem $4$ of Avgustinovich, Kitaev, and Taranenko \cite{akt}, a minimal $(k, \ell)$-right-crucial permutation has length $k+\ell-3$. Therefore, a next minimal $(k, \ell)$-right-crucial permutation has length $k+\ell-2$. By \Cref{right_crucial_char}, $P_{n}(\pi)$ must have $k-1$ columns and $\ell-1$ rows. Therefore, the only possible shape of $P_{n}(\pi)$ is $(k-1, 2, \underbrace{1, \ldots, 1}_{\ell-3})$. 

We first fill $P_{n}(\pi)$. By \Cref{lem_for_next_2}, there are
\begin{equation} \label{right_eq_1}
    k+\ell-4
\end{equation} ways to fill $P_{n}(\pi)$ such that it satisfies condition 3 of \Cref{right_crucial_char}.  

Now, we fill $Q_{n}(\pi)$. By the hook-length formula, there are 
\begin{equation} \label{right_eq_2}
    \frac{(k-2)(\ell-2)}{k+\ell-3} {\ell+k-2 \choose \ell-1}
\end{equation} ways to fill $Q_{n}(\pi)$. 

By the RSK correspondence, we finish by multiplying \Cref{right_eq_1} and \Cref{right_eq_2}.
\end{proof} 

%\section{Examples}
%\subfile{sections/examples.tex}

\section{Bicrucial Permutations}
\label{bi}
In this section, we enumerate the next minimal $(k, \ell)$-bicrucial permutations and the $(k, 3)$-bicrucial permutations. We begin by citing the analog of \Cref{right_crucial_char} for $(k, \ell)$-bicrucial permutations by Avgustinovich, Kitaev, and Taranenko \cite{akt}. 

\begin{thm}[Avgustinovich, Kitaev, and Taranenko \cite{akt}] \label{bicrucial_char}
The permutation $\pi \in S_{n}$ is $(k,\ell)$-bicrucial if and only if:
\begin{enumerate}
    \item the number of columns in $P_{n}(\pi)$ (and $Q_{n}(\pi)$) is $k-1$;
    \item the number of rows in $P_{n}(\pi)$ (and $Q_{n}(\pi)$) is $\ell-1$;
    \item $P_{n}(\pi)$ contains an increasing sequence of elements $p_{1, k-1}, p_{2, j_2}, \ldots, p_{\ell-1, j_{\ell-1}}$ for some $j_2,$ $\ldots,$ $j_{\ell-1}$; 
    \item $P_{n}(\pi)$ contains an increasing sequence of elements $p_{\ell-1,1}, p_{i_2, 2}, \ldots, p_{i_{k-1}, k-1}$ for some $i_2,$ $\ldots,$ $i_{k-1}$. 
\end{enumerate}
\end{thm}

\subsection{Next minimal bicrucial permutations}

To count the next minimal $(k, \ell)$-bicrucial permutations, we first establish the following lemma that counts the distinct $P_{n}(\pi)$ that satisfy $p_{\ell-1, 1} < p_{1, k-1}$ for some next minimal $(k, \ell)$-bicrucial permutation $\pi$. 

\begin{lem} \label{lem_for_next_3}
The number of standard Young tableaux $T$ of shape $(k-1, 3, \underbrace{2, \ldots, 2}_{\ell-3})$ that satisfy $t_{1,k-1} < t_{2, j} < t_{3,2}$ for some $j$ and $t_{\ell-1,1} < t_{i_2, 2} < t_{i_3,3} < t_{1, 4}$ for some $i_2$ and $i_3$ is $k + 2\ell-7$. 
\end{lem}
\begin{proof}
We begin by showing that $t_{1,k-1} \in \{k+\ell-3, k+\ell-2\}$. We first observe that because $t_{1, k-1} < t_{2, j} < t_{3,2}$, the entry $t_{1, k-1}$ must be less than the entries $\{t_{2, j}\} \cup \{t_{3,2},t_{4,2}, \ldots, t_{\ell-1, 2}\}$. Since at least $\ell-2$ entries are greater than $t_{1, k-1}$, the entry $t_{1, k-1}$ must be at most $k+\ell-2$. At the same time, because $t_{\ell-1,1} < t_{1, 4}$, the entry $t_{1, k-1}$ must be greater than the entries $\{t_{1,1}, t_{1, 2}, \ldots, t_{k-2}\} \cup \{t_{1,1}, t_{2,1}, \ldots, t_{\ell-1,1}\}$. Therefore, at least $k+\ell-4$ entries are less than $t_{1, k-1}$. So, the entry $t_{1, k-1}$ must be at least $k+\ell-3$. Furthermore, if $t_{1, k-1} = k+\ell-3$, then $\{t_{1,1}, t_{1, 2}, \ldots, t_{k-2}\} \cup \{t_{1,1}, t_{2,1}, \ldots, t_{\ell-1,1}\} = \{1, 2, \ldots, k+\ell-4\}$. 

 As we proceed, let $S$ be the filled Young diagram formed by the boxes $\{b_{1, 1}, b_{1, 2}, \ldots, b_{1, k-1}\} \cup \{b_{1,1}, b_{2, 1}, \ldots, b_{\ell-1, 1}\}$. Refer to \Cref{s_in_t_again} for the position of $S$ in $T$.
 
 \begin{figure}
\begin{tikzpicture}[scale=0.4]

\draw [line width=0.5](0,0)--(2,0);
\draw [line width=0.5](0,0)--(0,8);
\draw [line width=0.5](0,8)--(10,8);
\draw [line width=0.5](10,8)--(10,7);
\draw [line width=0.5](10,7)--(1,7);
\draw [line width=0.5](2,0)--(2,6);
\draw [line width=0.5](1,0)--(1,7);
\draw [line width=0.5](2,6)--(3,6);
\draw [line width=0.5](3,6)--(3,7);

\node [below] at (-2,5){{\large T = }};
\node [below] at (2,7.2){{\small $T \setminus S$}};

      %% brace a top
      
\node [below] at (0.5,8){{\small S}};
\end{tikzpicture} 
\caption{The position of $S$ in $T$}
\label{s_in_t_again}
\end{figure}
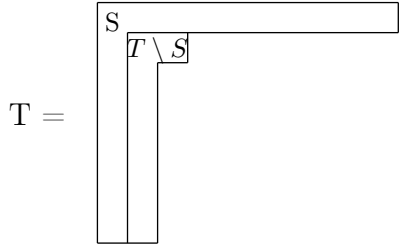
 
 First, consider the case $t_{1, k-1} = k+\ell-3$. As deduced previously, we have that $\{t_{1,1}, t_{1, 2}, \ldots, t_{k-2}\}$ $\cup$ $\{t_{1,1}, t_{2,1}, \ldots, t_{\ell-1,1}\} = \{1, 2, \ldots, k+\ell-4\}$. Thus, any entry in $T \setminus S$ must be greater than $t_{1, k-1}$. Therefore, the partial orders between an entry in $S$ and an entry in $T \setminus S$ given by the shape of $T$ are always satisfied. Thus, we fill $S$ and $T \setminus S$ independently.
 
 We begin by filling $S$. Because any entry in $T \setminus S$ is greater than any entry in $S$, for $t_{\ell-1, 1} < t_{i_2, 2} < t_{i_3, 3} < t_{1, 4}$ for some $i_2$ and $i_3$, it must be that $i_2 = i_3 = 1$. Because this implies that
 \begin{equation*}
     t_{1,1} < t_{2, 1} < \cdots < t_{\ell-1, 1} < t_{1, 2} < t_{1, 3} < \cdots < t_{1, k-1},
 \end{equation*} all entries of $S$ must be fixed. 
 
 We now fill $T \setminus S$. Because $t_{1, k-1} < t_{2,2} < t_{3,2}$, the only requirement when filling $T \setminus S$ is that $T$ is a standard Young tableau. Therefore, by the hook-length formula, there are 
\begin{equation} \label{lem_7_eq_1}
    \ell-2
\end{equation} ways of filling $T \setminus S$. 

Now, consider the case $t_{1, k-1} = k+\ell-2$. From $t_{\ell-1} < t_{1, 4}$, every entry of $S$ must still be less than $t_{1, k-1}$. But now, because $t_{1, k-1} = k+\ell-2$, and there are $k+\ell-3$ boxes in $S$, exactly one entry of $T \setminus S$ must be less than $t_{1, k-1}$. For $T$ to be a standard Young tableau, $t_{2,2}$ must be the minimum entry of $T \setminus S$. Thus, $t_{2,2}$ must be the entry in $T \setminus S$ less than $t_{1, k-1}$. Now, because it must be that $t_{1,k-1} < t_{2, 3} < t_{3, 2}$, we have that $t_{2, 3} < t_{3, 2} < t_{4,2} < \cdots < t_{\ell-1, 2}$. Therefore, the entries of $T \setminus (S \cup \{b_{2,2}\})$ must be fixed. Now, let $U$ be the standard Young tableau formed by the boxes $(S \cup \{b_{2,2}\}) \setminus \{b_{1, k-1}\}$. Refer to \Cref{u_in_t_again} for the position of $U$ in $T$. It now remains to fill $U$ with $\{1, 2, \ldots, \ell+k-3\}$ such that $t_{\ell-1,1} < t_{i_2,2} < t_{1,3}$ for some $i_2$. Therefore, by \Cref{lem_for_next_3}, there are 
\begin{equation} \label{lem_7_eq_2}
    k+\ell-5
\end{equation} ways of filling $U$. 

We finish by summing \Cref{lem_7_eq_1} and \Cref{lem_7_eq_2}. 
 \begin{figure}[h]
\begin{tikzpicture}[scale=0.4]

\draw [line width=0.5](0,0)--(2,0);
\draw [line width=0.5](0,0)--(0,8);
\draw [line width=0.5](0,8)--(10,8);
\draw [line width=0.5](10,8)--(10,7);
\draw [line width=0.5](9,8)--(9,7);
\draw [line width=0.5](10,7)--(2,7);
\draw [line width=0.5](2,0)--(2,7);
\draw [line width=0.5](1,0)--(1,6);
\draw [line width=0.5](1,6)--(3,6);
\draw [line width=0.5](3,6)--(3,7);

\node [below] at (-2,5){{\large T = }};

      %% brace a top
      
\node [below] at (0.5,8){{\small U}};
\end{tikzpicture} 
\caption{The position of $U$ in $T$}
\label{u_in_t_again}
\end{figure}
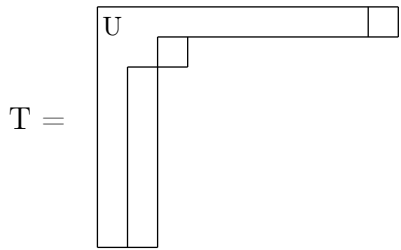
\end{proof} 

Next, we take the transpose of the standard Young tableaux in \Cref{lem_for_next_3} to arrive at the following lemma. 

\begin{lem} \label{lem_next_transpose}
For $k, \ell \geq 3$, the number of standard Young tableaux $T$ of shape $(k-1, k-1, 2, \underbrace{1, \ldots, 1}_{\ell-4})$ that satisfy $t_{1,k-1} < t_{2, j_2} < t_{3,j_3} < t_{4,1}$ for some $j_2$ and $j_3$ and $t_{\ell-1,1} < t_{i, 2} < t_{2, 3}$ for some $i$ is $2k + \ell-7$. 
\end{lem}
\begin{proof}
Take the transpose of $T$ to arrive at $T^{t}$. Refer to \Cref{p_shape} for the shape of $T$ and $T^{t}$. 

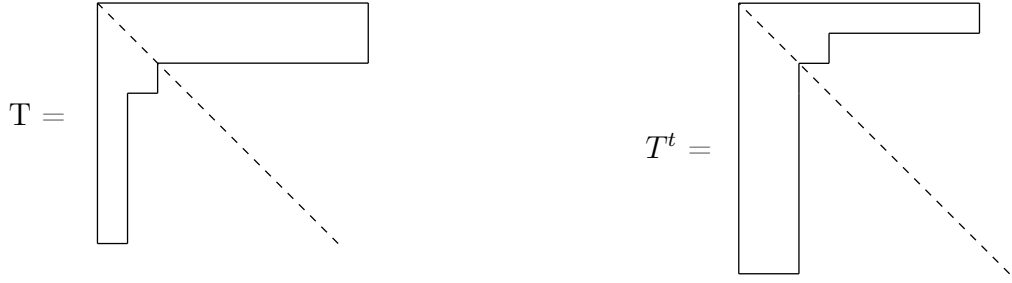
\begin{figure}[h]
\begin{multicols}{2}
\center 
\begin{tikzpicture}[scale=0.4]

\draw [line width=0.5](0,0)--(1,0);
\draw [line width=0.5](0,0)--(0,8);
\draw [line width=0.5](0,8)--(9,8);
\draw [line width=0.5](9,8)--(9,6);
\draw [line width=0.5](1,0)--(1,5);
\draw [line width=0.5](2,5)--(2,6);
\draw [line width=0.5](2,5)--(1,5);
\draw [line width=0.5](2,6)--(9,6);
\draw [line width=0.5, dashed](8,0)--(0,8);

\node [below] at (-2,5){{\large T = }};

      %% brace a top
      
\end{tikzpicture} \columnbreak
\center 
\begin{tikzpicture}[scale=0.4]

\draw [line width=0.5](0,0)--(2,0);
\draw [line width=0.5](0,0)--(0,9);
\draw [line width=0.5](0,9)--(8,9);
\draw [line width=0.5](8,9)--(8,8);
\draw [line width=0.5](2,0)--(2,6);
\draw [line width=0.5](2,6)--(2,7);
\draw [line width=0.5](2,7)--(3,7);
\draw [line width=0.5](3,8)--(3,7);
\draw [line width=0.5](3,8)--(8,8);
\draw [line width=0.5, dashed](9,0)--(0,9);

\node [below] at (-2,5){{\large $T^{t}$ = }};

      %% brace a top
      
\end{tikzpicture}
\end{multicols} 
\caption{The construction of $T^{t}$ from $T$} \label{p_shape}
\end{figure}

Because the transformation from $T$ to $T^{t}$ simply switches the indices of rows and columns, $T^{t}$ must be a standard Young tableau that satisfies the conditions of \Cref{lem_for_next_3}. Therefore, there must be $2k + \ell -7$ ways of filling $T$. 
\end{proof}

Finally, we establish the following lemma that counts the number of distinct $P_{n}(\pi)$ that satisfy $p_{1, k-1}< p_{\ell-1, 1}$ for some next minimal $(k, \ell)$-bicrucial permutation $\pi$ when $k = \ell+1$. 

\begin{lem} \label{lem_next_one}
The number of standard Young tableaux of shape $(k-1, k-1, \underbrace{1, \ldots, 1}_{k-4})$ that satisfy $t_{1, k-1} < t_{2,j} < t_{3, 1}$ for some $j$ and $t_{k-2, 1} < t_{2,2}$ is $1$. 
\end{lem}
\begin{proof}
We begin by letting $S$ be the filled Young diagram formed by the boxes $\{b_{1,1}, b_{1, 2}, \ldots, b_{1, k-1}\}$ and $\{b_{1,1}, b_{2, 1}, \ldots, b_{\ell-1, 1}\}$. Refer to \Cref{fig:r_in_p} for the position of $S$ in $T$.  
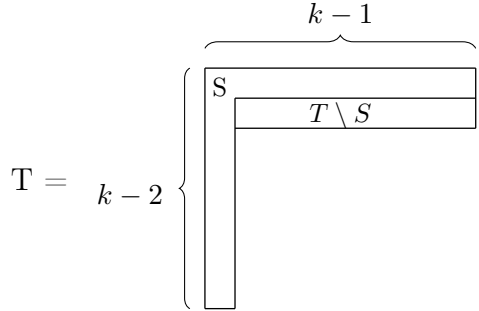
\begin{figure}[h]
    \centering
    \begin{tikzpicture}[scale=0.4]

\draw [line width=0.5](0,0)--(1,0);
\draw [line width=0.5](0,0)--(0,8);
\draw [line width=0.5](0,8)--(9,8);
\draw [line width=0.5](9,8)--(9,6);
\draw [line width=0.5](1,0)--(1,7);
\draw [line width=0.5](1,6)--(9,6);
\draw [line width=0.5](1,7)--(9,7);
\node [below] at (-5.5,5){{\large T = }};
\node [below] at (0.5,8){{\small S}};
\node [below] at (4.5,7.2){{\small $T \setminus S$}};

\draw [decorate,decoration={brace,amplitude=4pt, mirror},xshift=0.5cm,yshift=0pt]
      (-1,8) -- (-1,0) node [midway,right,xshift=.1cm] {};      
\node [below] at (-2.5,4.5){{$k-2$}};  

\draw [decorate,decoration={brace,amplitude=5pt,raise=4ex}]
  (0,7) -- (9,7) node[midway,yshift=-3em]{};
\node [below] at (4.5,10.5){{$k-1$}}; 
      %% brace a top
      
\end{tikzpicture}
    \caption{The position of $S$ in $T$}
    \label{fig:r_in_p}
\end{figure} 

First, because $t_{1, k-1} < t_{3, 1}$ from the statement of the lemma and $t_{3, 1}< t_{k-2, 1}$ from $T$ being a standard Young tableau, we have that $t_{1, k-1} < t_{k-2,1}$. Thus, $t_{k-2, 1}$ must be the maximum entry of $S$. At the same time, from $t_{k-2, 1} < t_{2,2}$, we have that each entry of $T \setminus S$ must be greater than $t_{k-2, 1}$. Therefore, to satisfy $t_{1, k-1} < t_{j_2, 2} < t_{3, 1}$ for some $j_2$, we must have $t_{1, k-1} < t_{2, 1}$. Now, collating $t_{1, k-1} < t_{2, 1}$ and $t_{k-2, 1} < t_{2,2}$, we have a complete ordering on the entries of $T$:
\begin{equation*}
    t_{1,1} < t_{1, 2} < \cdots < t_{1, k-1} < t_{2, 1} < \cdots < t_{k-2, 1} < t_{2,2} < t_{2, 3} \cdots < t_{2, k-1}. 
\end{equation*}
Therefore, the filling of $T$ is fixed as sought.
\end{proof}

Our enumeration of next minimal $(k, \ell)$-bricrucial permutations now follows from the hook-length formula, \Cref{lem_for_next_2}, \Cref{bicrucial_char}, \Cref{lem_for_next_3}, \Cref{lem_next_transpose}, and \Cref{lem_next_one}. We deal with the cases $(k, \ell) \in \{(4, 3), (4,4)\}$ separately, because unlike $(k, \ell)$ with larger $k$ for which $k = \ell$ or $k = \ell+1$, the shape of $P_{n}(\pi)$ (and $Q_{n}(\pi)$) are fixed when $(k, \ell) \in \{(4, 3), (4,4)\}$. 

\begin{thm} \label{bicrucial_next}
Suppose that $k \geq \ell$. The number of next minimal $(4,3)$-bicrucial permutations is $25$, the number of next minimal $(4,4)$-bicrucial permutations is $756$, and the number of $(k, \ell)$-bicrucial permutations for $(k, \ell) \not \in \{(4,3), (4,4)\}$ is 
\begin{footnotesize}
\begin{align*}
    (1 + \delta_{k, \ell}) \frac{(k+2\ell-7)(k+2\ell-4)(k-3)}{(k+\ell-4)(k+\ell-3)} \binom{k+2\ell-5}{k-2, \ell-3, \ell} + \delta_{k-1, \ell} \frac{3k-6}{(2k-5)(2k-4)} \binom{3k-7}{k-4, k-2, k-1}.
\end{align*}
\end{footnotesize}
\end{thm}
\begin{proof}
We begin by supposing that $(k, \ell) \not \in \{(4,3), (4,4)\}$. By Theorem 20 of Avgustinovich, Kitaev, and Taranenko \cite{akt}, the length of a minimal $(k, \ell)$-bicrucial permutation is $k+2 \ell-5$. Therefore, the length of a next minimal $(k, \ell)$-bicrucial permutations is $k + 2 \ell-4$. In addition, by conditions $1$ and $2$ of \Cref{bicrucial_char}, $P_{n}(\pi)$ must have $k-1$ columns and $\ell-1$ rows. 

First, suppose that $p_{\ell-1, 1} < p_{1, k-1}$. Because all entries in column $1$ are less than $p_{1, k-1}$, to satisfy condition 3 of \Cref{bicrucial_char}, column $2$ of $P_{n}(\pi)$ must have $\ell-1$ boxes. Thus, for $P_{n}(\pi)$ to have $k-1$ boxes in row $1$, $\ell-1$ boxes in columns $1$ and $2$, and $k+2\ell-4$ boxes total, $P_{n}(\pi)$ must have the shape $(k-1, 3, \underbrace{2, \ldots, 2}_{\ell-3})$. 

We first fill $P_{n}(\pi)$. First, by \Cref{lem_for_next_3}, there are 
\begin{equation} \label{next_min_bi_eq_4}
    k + 2 \ell - 7
\end{equation} ways to fill $P_{n}(\pi)$ such that it satisfies conditions $3$ and $4$ of \Cref{bicrucial_char}. 

Now, we fill $Q_{n}(\pi)$. The only restriction on $Q_{n}(\pi)$ is that it is a standard Young tableau that has the same shape as $P_{n}(\pi)$. Therefore, by the hook-length formula, there are 
\begin{equation} \label{next_min_bi_eq_5}
    \frac{(k+2\ell-4)(k-3)}{(k+\ell-4)(k+\ell-3)} \binom{k+2\ell-5}{k-2, \ell-3, \ell}
\end{equation}
ways to fill $Q_{n}(\pi)$.

By the RSK correspondence, we finish by multiplying \Cref{next_min_bi_eq_4} and \Cref{next_min_bi_eq_5}.

Now, suppose that $p_{1, k-1} < p_{\ell-1, 1}$. Then, because all entries in row $1$ are less than $p_{\ell-1, 1}$, for condition $4$ of \Cref{bicrucial_char} to hold, row $2$ must have all $k-1$ boxes. Therefore, the minimum number of boxes in $P_{n}(\pi)$ is $2k + \ell-5$. Since $\pi$ has length $k + 2 \ell- 4$, it must be that $k \in \{\ell, \ell+1\}$. 

Suppose that $k = \ell$. Then, for $P_{n}(\pi)$ to have $k-1$ boxes in rows $1$ and $2$, $k-1$ boxes in column $1$, and $3k-4$ boxes total, $P_{n}(\pi)$ must have the shape $(k-1, k-1, 2, \underbrace{1, \ldots, 1}_{k-4})$. Therefore, by \Cref{lem_next_transpose}, there are
\begin{equation} \label{next_min_bi_eq_1}
    3k-7    
\end{equation} ways of filling $P_{n}(\pi)$. Now, the only restriction on $Q_{n}(\pi)$ is that it is a standard Young tableau of the same shape as $P_{n}(\pi)$. Therefore, by the hook-length formula, there are
\begin{equation} \label{next_min_bi_eq_2}
    \frac{(3k-4)(k-3)}{(2k-4)(2k-3)} \binom{k+2\ell-5}{k-3, k-2, k}
\end{equation} ways to fill $Q_{n}(\pi)$. 

By the RSK correspondence, we finish this case by multiplying \Cref{next_min_bi_eq_1} and \Cref{next_min_bi_eq_2}. 

Now, if $k = \ell+1$, then for $P_{n}(\pi)$ to have $k-1$ boxes in rows $1$ and $2$, $\ell-1$ boxes in column $1$, and $k+2\ell-4$ boxes total, $P_{n}(\pi)$ must have shape $(k-1, k-1, \underbrace{1, \ldots, 1}_{\ell-3})$. Therefore, by \Cref{lem_next_one}, there is exactly one way to fill $P_{n}(\pi)$. 

Next, we fill $Q_{n}(\pi)$. The only restriction on $Q_{n}(\pi)$ is that it is a standard Young tableau of the same shape as $P_{n}(\pi)$. Therefore, by the hook-length formula, there are 
\begin{equation} \label{next_min_bi_eq_3}
    \frac{3k-6}{(2k-5)(2k-4)} \binom{3k-7}{k-4, k-2, k-1} 
\end{equation} ways of filling $Q_{n}(\pi)$. 

By the RSK correspondence, we see that \Cref{next_min_bi_eq_3} gives the number of next minimal $(k, k-1)$-tricrucial permutations in this case. 

To complete the enumeration of the next minimal $(k, \ell)$-bicrucial permutations, we ran a quick computer search to show that the number of next minimal $(4,3)$-bicrucial permutations is $25$ and the number of next minimal $(4,4)$-bicrucial permutations is $756$. 
\end{proof}

\subsection{$(k, 3)$-bicrucial permutations}

To count the $(k, 3)$-bicrucial permutations, we first establish the following lemma that counts the number of distinct $P_{n}(\pi)$ across $(k, 3)$-bicrucial permutations $\pi$. 

\begin{lem} \label{lem_for_min}
The number of standard Young tableaux $T$ of shape $(k-1, \underbrace{2, \ldots, 2}_{\ell-2})$ for which the exact set of entries less than $t_{1, k-1}$ is $\{t_{1, 1}, t_{1, 2}, \ldots, t_{1, k-2}\} \cup \{t_{1,1}, \ldots, t_{i, 1}\}$ is $$\frac{i}{2\ell-i-2} {k+i-4 \choose i-1}{2 \ell-i-2 \choose \ell-1}.$$
\end{lem}
\begin{proof}
Let $S$ be the standard Young tableau formed by the boxes $\{b_{1, 1}, b_{1, 2}, \ldots, b_{1, k-2}, b_{1, k-1}\} \cup \{b_{1,1}, \ldots, b_{i, 1}\}$. Refer to \Cref{fig:t_t''} for the position of $S$ in $T$. 

Since the elements in $T\setminus S$ are always greater than that of $S$, the partial orders between an element in $S$ and an element in $T \setminus S$ are always respected. It thus remains to fill $S$ and $T \setminus S$ independently such that any partial order given by the shape of $T$ between two elements of $S$ and between two elements of $T \setminus S$ are respected. 
\begin{figure}[h]
\center 
\begin{tikzpicture}[scale=0.4]

\draw [line width=0.5](0,0)--(2,0);
\draw [line width=0.5](0,0)--(0,8);
\draw [line width=0.5](0,8)--(10,8);
\draw [line width=0.5](10,8)--(10,7);
\draw [line width=0.5](10,7)--(1,7);
\draw [line width=0.5](2,0)--(2,7);
\draw [line width=0.5](1,7)--(1,5);
\draw [line width=0.5](0,5)--(1,5);

\node [below] at (-6,5){{\large T = }};

 \draw [decorate,decoration={brace,amplitude=4pt, mirror},xshift=0.5cm,yshift=0pt]
      (-1,8) -- (-1,5) node [midway,right,xshift=.1cm] {};

 \draw [decorate,decoration={brace,amplitude=4pt, mirror},xshift=0.5cm,yshift=0pt]
      (-1,5) -- (-1,0) node [midway,right,xshift=.1cm] {};      
\node [below] at (-3,3){{$\ell-i-1$}};      
      %% brace a top
      
\node [below] at (0.5,8){{\small S}};
\node [below] at (-1.5,6.5){{$i$}};
\node [below] at (1,3.5){{\small $T \setminus S$}};

\end{tikzpicture} 
\caption{The position of $S$ in $T$} \label{fig:t_t''}
\end{figure}
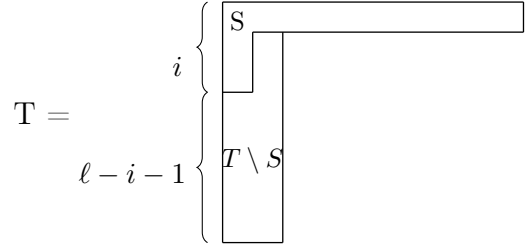

We first fill $S$. We have that $t_{1, k-1} = k+i-2$, and $\{t_{1, 1}, t_{1, 2}, \ldots, t_{1, k-2}\} \cup \{t_{1,1}, \ldots, t_{i, 1}\} = \{1, 2, \ldots, k+i-3\}$. Thus, by the hook-length formula, there are\begin{equation} \label{lem_2_eq_1}
    {k+i-4 \choose i-1}
\end{equation} ways of filling $S$. 

Next, we fill $T\setminus S$. To do so, we rotate the complement of $T \setminus S$ by $180$ degrees. Then, we have a filled Young diagram of shape $(\underbrace{2, 2, \ldots, 2}_{\ell-i-1}, \underbrace{1, 1, \ldots, 1}_{i-1})$ as each row and column must be increasing. Therefore, by the hook-length formula, there are \begin{equation} \label{lem_2_eq_2}
    \frac{i}{2\ell-i-2} {2\ell-i-2 \choose \ell-1}
\end{equation} ways of filling $T \setminus S$. 

We finish by multiplying \Cref{lem_2_eq_1} and \Cref{lem_2_eq_2}. 
\end{proof}

Our enumeration of $(k, 3)$-bicrucial permutations now follows from the hook-length formula, \Cref{bicrucial_char}, and \Cref{lem_for_min}.

\begin{thm} \label{bicrucial_3}
Let $k \geq 3$. The number of $(k, 3)$-bicrucial permutations of length $n$ for $k+1 \leq n \leq 2k-3$ is
$$\frac{2k-n-1}{n+1} {n+1 \choose k} \left( \frac{2k-n}{n}{n \choose k} - \frac{2k-n-1}{n-1} {n-1 \choose k-1} \right).$$ For $n = 2(k-1)$, it is $$\left( \frac{1}{k} {2k-2 \choose k-1} \right)^2.$$ 
\end{thm}
\begin{proof}
Since $\ell = 3$, the third condition of \Cref{bicrucial_char} is equivalent to $p_{1, k-1} < p_{2, n-k+1}$. Since the rows of $P_{n}(\pi)$ increase, either $p_{1, k-1} = n$ or $p_{2, n-k+1} = n$. Therefore, the third condition of \Cref{bicrucial_char} is equivalent to $p_{2, n-k+1}=n$. 

Now, suppose that $k+1 \leq n \leq 2k-3$. Then by \Cref{lem_for_min}, there are
\begin{equation} \label{bicrucial_eq_1}
  \frac{2k-n}{n} {n \choose k} - \frac{2k-n-1}{n-1} {n-1 \choose k-1}  
\end{equation} possible fillings of $P_{n}(\pi)$. 

We now fill $Q_{n}(\pi)$. From \Cref{bicrucial_char}, the only restriction on $Q_{n}(\pi)$ is that it is a standard Young tableau of shape $(k-1, n-k+1)$. Therefore, by the hook-length formula, there are \begin{equation} \label{bicrucia_eq_2}
    \frac{2k-n-1}{n+1} {n+1 \choose k}
\end{equation} possible fillings of $Q_{n}(\pi)$. 

By the RSK correspondence, we finish by multiplying \Cref{bicrucial_eq_1} and \Cref{bicrucia_eq_2}. 

Now, suppose that $n = 2(k-1)$. Then the third condition of \Cref{bicrucial_char} is always met, because $p_{1, k-1} < p_{2, k-1} = n$. The fourth condition of \Cref{bicrucial_char} is also always met, because $p_{2, 1} < p_{2, 2} < \cdots < p_{2, k-1}$. Thus, the only restrictions on $P_{n}(\pi)$ and $Q_{n}(\pi)$ are that they are standard Young tableaux of shape $(k-1, k-1)$. Therefore, by the hook-length formula, there are \begin{equation} \label{bicrucial_eq_3}
    \frac{1}{k}{2k-2 \choose k-1}
\end{equation} possible fillings for each of $P_{n}(\pi)$ and $Q_{n}(\pi)$. 

By the RSK correspondence, we finish by taking the square of \Cref{bicrucial_eq_3}. 
\end{proof}

That all $(k, \ell)$-bicrucial permutations for $k \geq \ell$ must have length at least $k+2\ell-5$ and at most $(k-1)(\ell-1)-1$ follows from Theorem 9 of Avgustinovich, Kitaev, and Taranenko \cite{akt}. Therefore, \Cref{bicrucial_3} gives a complete enumeration of $(k, 3)$-bicrucial permutations. 

\section{Tricrucial Permutations}
\label{tri}
In this section, we enumerate the minimal $(k, \ell)$-tricrucial permutations, the next minimal $(k, \ell)$-tricrucial permutations, and the $(k, 3)$-tricrucial permutations. We begin by citing the analog of \Cref{right_crucial_char} for $(k, \ell)$-tricrucial permutations. 

\begin{thm}[Avgustinovich, Kitaev, and Taranenko \cite{akt}] \label{tricrucial_char}
The permutation $\pi$ is $(k,\ell)$-tricrucial if and only if:
\begin{enumerate}
    \item the number of columns in $P_{n}(\pi)$ (and $Q_{n}(\pi)$) is $k-1$;
    \item the number of rows in $P_{n}(\pi)$ (and $Q_{n}(\pi)$) is $\ell-1$;
    \item $P_{n}(\pi)$ contains an increasing sequence of elements $p_{1, k-1}, p_{2, j_2}, \ldots, p_{\ell-1, j_{\ell-1}}$ for some $j_2, \ldots,$ $j_{\ell-1}$; 
    \item $P_{n}(\pi)$ contains an increasing sequence of elements $p_{\ell-1,1}, p_{i_2, 2}, \ldots,p_{i_{k-1}, k-1}$ for some $i_2, \ldots,$ $i_{k-1}$;
    \item $Q_{n}(\pi)$ contains an increasing sequence of elements $q_{1, k-1}, q_{2, m_2}, \ldots, q_{\ell-1, m_{\ell-1}}$ for some $m_2, $ $\ldots,$ $m_{\ell-1}$. 
\end{enumerate}
\end{thm}

\subsection{Minimal tricrucial permutations}

When $\pi$ is minimal $(k, \ell)$-tricrucial, Avgustinovich, Kitaev, and Taranenko \cite{akt} showed that the characterization of $(P_{n}(\pi), Q_{n}(\pi))$ presented in \Cref{tricrucial_char} is equivalent to the following. 

\begin{thm}[Avgustinovich, Kitaev, and Taranenko \cite{akt}] \label{tricrucial_rsk}
Let $T$ be the standard Young tableau of shape $(k-1, \underbrace{2, \ldots, 2}_{\ell-2})$,
\begin{itemize}
    \item $t_{1, i} = i$ for $1 \leq i \leq \ell-1$;
    \item $t_{2, i} = i + k + \ell -4$ for $2 \leq i \leq \ell-1$; 
    \item $t_{i, 1} = \ell-1 + i$ for $2 \leq i \leq k-1$;
\end{itemize}. 
\begin{center}
\begin{tikzpicture}[scale=0.45]

\draw [line width=0.5](0,5)--(7.3,5);
\draw [line width=0.5](8.7,5)--(12,5);
\draw [line width=0.5](0,4)--(7.3,4);
\draw [line width=0.5](8.7,4)--(12,4);
\draw [line width=0.5](0,3)--(5,3);
\draw [line width=0.5](0,1)--(5,1);
\draw [line width=0.5](0,0)--(5,0);

\draw [line width=0.5](0,2.7)--(0,5);
\draw [line width=0.5](0,0)--(0,1.3);
\draw [line width=0.5](1.6,2.7)--(1.6,5);
\draw [line width=0.5](1.6,0)--(1.6,1.3);
\draw [line width=0.5](5,0)--(5,1.3);
\draw [line width=0.5](5,2.7)--(5,5);
\draw [line width=0.5](7,4)--(7,5);
\draw [line width=0.5](9,4)--(9,5);
\draw [line width=0.5](12,4)--(12,5);

\node [below] at (-3,3){{\large T = }};
\node [below] at (0.8,5){{\small $1$}};
\node [below] at (3,5){{\small $\ell$}};
\node [below] at (6,5){{\small $\ell+1$}};
\node [below] at (8,5){{\small $\cdots$}};
\node [below] at (10.5,5){{\small $k+\ell-3$}};
\node [below] at (3.1,4){{\small $k+\ell-2$}};
\node [below] at (0.8,4){{\small $2$}};
\node [below] at (0.8,3.2){{\small $\vdots$}};
\node [below] at (3.2,3.2){{\small $\vdots$}};
\node [below] at (0.8,1){{\small $\ell-1$}};
\node [below] at (3.3,1){{\small $k+2\ell-5$}};
\end{tikzpicture}
\end{center}
When $k > \ell$, $\pi$ is $(k, \ell)$-tricrucial if and only if $P_{n}(\pi)$ is given by $T$, and $Q_{n}(\pi)$ of the same shape satisfies $q_{1, k-1} < q_{2, 2}$. 

When $k = \ell$, $\pi$ is $(k, \ell)$-tricrucial if and only if $P_{n}(\pi)$ is given by $T$
or its transpose, and $Q_{n}(\pi)$ of the same shape satisfies $q_{1, k-1} < q_{2, 2}$ in the former case, and $q_{1, k-1} < q_{2, j} < q_{3,1}$ for some $j$ in the latter case. 
\end{thm}

We next establish the following lemma, which counts the number of $Q_{n}(\pi)$ that has shape $(k-1, \underbrace{2, \ldots, 2}_{k-2})$ and $\pi$ is some minimal $(k, \ell)$-tricrucial permutation. 

\begin{lem} \label{bijection}
There exists an explicit bijection between the set of standard Young tableaux $T$ of shape $(k-1, \underbrace{2, \ldots, 2}_{\ell-2})$ that satisfy $t_{1,k-1} < t_{2,2}$ and the set of standard Young tableaux $T'$ of shape $(\ell-1, \ell-1, \underbrace{1, \ldots, 1}_{k-3})$ that satisfy $t'_{1, \ell-1} < t'_{2, j} < t'_{3,1}$ for some $j$.
\end{lem}
\begin{proof}
Because $T$ must satisfy $t_{1, k-1} < t_{2, 2}$, the exact set of entries that are less than $t_{1, k-1}$ must be $\{t_{1, 1}, t_{2, 1}, \ldots, t_{i, 1}\} \cup \{t_{1, 2}, t_{1, 3}, \ldots, t_{1, k-2}\}$ for some $i$. Similarly, as $T'$ must satisfy $t'_{1, \ell-1} < t'_{2, j} < t'_{3,1}$ for some $j$, the exact set of entries greater than $t'_{1, \ell-1}$ must be $\{t'_{2, \ell-i}, t'_{2, \ell-i+1}, \ldots, t'_{2, \ell-1}\} \cup \{t'_{3,1}, t'_{4,1}, \ldots, t'_{k-1, 1}\}$ for some $i$, and it must be that $t'_{2, \ell-i} = t'_{1, \ell-1} + 1$. 

We break down $T$ as follows such that the filled Young diagrams $B$ and $C$ contain the exact set of entries that are less than $t_{1, k-1}$: 
\begin{itemize}
    \item $A:$ the filled Young diagram formed by the box $\{b_{1, k-1}\}$; 
    \item $B:$ the filled Young diagram formed by the boxes $\{b_{1, 2}, b_{1,3}, \ldots, b_{1, k-2}\}$;
    \item $C:$ the filled Young diagram formed by the boxes $\{b_{1,1}, b_{2,1}, \ldots, b_{i, 1}\}$;
    \item $D:$ the filled skew Young diagram formed by the boxes $\{b_{i+1, 1}, b_{i+2, 1}, \ldots, b_{\ell-1, 1}\}$ and $\{b_{2,2}, b_{3,2}, \ldots, b_{\ell-1,2}\}$. 
\end{itemize}

Similarly, we break down $T'$ as follows such that the filled Young diagrams $B'$ and $C'$ contain the exact set of entries that are greater than $t'_{1, \ell-1}$:
\begin{itemize}
    \item $A':$ the filled Young diagram formed by the box $\{b'_{1, \ell-1}\}$;
    \item $B':$ the filled Young diagram formed by the boxes $\{b'_{3,1}, b'_{4,1}, \ldots, b'_{k-1, 1}\}$;
    \item $C':$ the filled Young diagram formed by the boxes $\{b'_{2, \ell-i}, b'_{2, \ell-i+1}, \ldots, b'_{2, \ell-1}\}$; 
    \item $D':$ the filled (skew) Young diagram formed by the boxes $\{b'_{1,1}, b'_{1,2}, \ldots, b'_{1, \ell-2}\}$ and \newline $\{b'_{2,1}, b'_{2,2}, \ldots, b'_{2, \ell-i-1}\}$. 
\end{itemize}

Refer to \Cref{fig:t_t'} for a visual breakdown of $T$ and $T'$. 

\begin{figure}[h]
\begin{multicols}{2}
\center 
\begin{tikzpicture}[scale=0.35]

\draw [line width=0.5](0,0)--(2,0);
\draw [line width=0.5](0,0)--(0,8);
\draw [line width=0.5](0,8)--(9,8);
\draw [line width=0.5](9,8)--(9,7);
\draw [line width=0.5](9,7)--(1,7);
\draw [line width=0.5](8,8)--(8,7);
\draw [line width=0.5](2,0)--(2,7);
\draw [line width=0.5](1,8)--(1,4);
\draw [line width=0.5](0,4)--(1,4);

\node [below] at (-6,5){{\large T = }};
 \draw [decorate,decoration={brace,amplitude=4pt, mirror},xshift=0.5cm,yshift=0pt]
      (-1,8) -- (-1,4) node [midway,right,xshift=.1cm] {};  

\node [below] at (-1.5,6.5){{$i$}};

 \draw [decorate,decoration={brace,amplitude=4pt},xshift=0.5cm,yshift=0pt]
      (-1,0) -- (-1,4) node [midway,right,xshift=.1cm] {};  

\draw [decorate,decoration={brace,amplitude=5pt,raise=4ex}]
  (0,6.6) -- (9,6.6) node[midway,yshift=-3em]{};

\node [below] at (4.5,10.5){{$k-1$}};

\node [below] at (-3,2.8){{$\ell-i-1$}};

\node [below] at (1,3){{\small D}};
\node [below] at (0.5,6.5){{\small C}};
\node [below] at (4,8.2){{\small B}};
\node [below] at (8.5,8.2){{\small A}};

\end{tikzpicture} \columnbreak
\center 
\begin{tikzpicture}[scale=0.35]

\draw [decorate,decoration={brace,amplitude=5pt,raise=4ex}]
  (0,7.6) -- (8,7.6) node[midway,yshift=-3em]{};

\node [below] at (4,11.5){{$\ell-1$}};

\draw [line width=0.5](0,0)--(1,0);
\draw [line width=0.5](0, 7)--(8,7);
\draw [line width=0.5](0, 0)--(0,7);
\draw [line width=0.5](1, 0)--(1,7);
\draw [line width=0.5](0, 7)--(0,9);
\draw [line width=0.5](0, 9)--(8,9);
\draw [line width=0.5](8, 9)--(8,7);
\draw [line width=0.5](4, 7)--(4,8);
\draw [line width=0.5](4, 8)--(8,8);
\draw [line width=0.5](7, 8)--(7,9);

 \draw [decorate,decoration={brace,amplitude=4pt, mirror},xshift=0.5cm,yshift=0pt]
      (-1,9) -- (-1,0) node [midway,right,xshift=.1cm] {};  

\node [below] at (-2.5,5.5){{$k-1$}};

\draw [decorate,decoration={brace,amplitude=5pt,mirror,raise=4ex}]
  (4,8.3) -- (8,8.3) node[midway,yshift=-3em]{$i$};

\node [below] at (-6,5){{\large T' = }};

\node [below] at (0.5,4.5){{\small B'}};
\node [below] at (2.5,8.7){{\small D'}};
\node [below] at (7.5,9.2){{\small A'}};
\node [below] at (6.2,8.2){{\small C'}};

\end{tikzpicture}
\end{multicols} 
\caption{Breakdown of $T$ and $T'$} \label{fig:t_t'}
\end{figure}
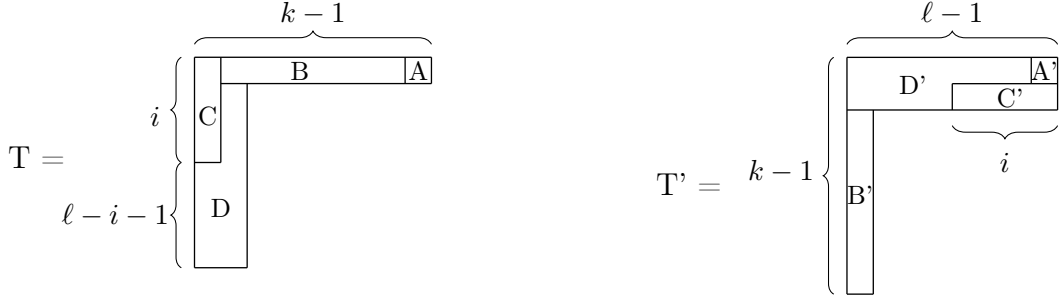

Now, consider the following map that maps a filled Young diagram $T'$ of shape $(\ell-1, \ell-1, \underbrace{1, \ldots, 1}_{k-3})$ that satisfies $t'_{2, \ell-i} = t'_{1, \ell-1}+1$ and contains the entries greater than $t'_{1, \ell-1}$ in $B'$ and $C'$ to a filled Young diagram $T$ of shape $(k-1, \underbrace{2, \ldots, 2}_{\ell-2})$ that contains the entries less than $t_{1, k-1}$ in $B$ and $C$: 
\begin{itemize}
    \item $t_{1, k-1} = (n+1)-t'_{1, \ell-1}$;
    \item $t_{r, 1} = t'_{2, r+\ell-i-1}-(2\ell-i-3)$ for $1 \leq r \leq i$;
    \item $t_{r, 1} = (n+1) - t'_{2, \ell-r}$ for $i+1 \leq r \leq \ell-1$; 
    \item $t_{r, 2} = (n+1) - t'_{1, \ell-r}$ for $2 \leq r \leq \ell-1$;
    \item $t_{1, r} = t'_{r+1, 1}-(2\ell-i-3)$ for $2 \leq r \leq k-2$. 
\end{itemize}

We can check that the proposed map maps $A'$ (respectively, $B'$, $C'$, and $D'$) to $A$ (respectively, $B$, $C$, and $D$). 

It remains to show that $T'$ is a standard Young tableau that satisfies $t'_{2, \ell-i} = t'_{1,\ell-1}+1$ and contains the entries greater than $t'_{1, \ell-1}$ in $B'$ and $C'$ if and only if $T$ is a standard Young tableau that contains the entries less than $t_{1, k-1}$ in $B$ and $C$. Note that because the entries of $B$ and $C$ are less than the entries of $A$ and $D$ in $T$, the partial orders between the entries of $A$ and $B$ (respectively, $A$ and $C$, $B$ and $D$, $C$ and $D$) are always satisfied. The same applies to the entries of $A'$ and $B'$ (respectively, $A'$ and $C'$, $B'$ and $D'$, $C'$ and $D'$) in $T'$.

Next, we first see that the partial orders between two elements that are both in $B$ (respectively, $C$, $D$) must be respected if and only if they are respected in $B'$ (respectively, $C'$, $D'$). Furthermore, $t'_{1, \ell-1} = 2\ell-i-3$, because in $T'$, the entries that are greater than $t'_{1, \ell-1}$ are contained in $B'$ and $C'$. Therefore, that $t'_{2, \ell-i} = t'_{1,\ell-1} + 1 = 2\ell -i-2$ is equivalent to $t_{1,1} = 1$ under the proposed map. But $t_{1,1}=1$ if and only if the only partial order given by the structure of $T$ between the elements of $B$ and $C$ are satisfied. Similarly, $t_{1, k-1} < t_{2,2}$ if and only if $t'_{1, \ell-2} = (n+1)-t_{2,2} < t'_{1, \ell-1} = (n+1)-t_{1, k-1}$. But $t'_{1, \ell-2} < t'_{1, \ell-1}$ if and only if the only partial order between the entries of $A'$ and $D'$ given by the structure of $T'$ is satisfied. Thus, collating, $T$ is a standard Young tableau that satisfies $t'_{2, \ell-i} = t'_{1,\ell-1}+1$ and contains the entries greater than $t'_{1, \ell-1}$ in $B'$ and $C'$ if and only if $T$ is a standard Young tableau that contains the entries less than $t_{1, k-1}$ in $B$ and $C$ as sought.  
\end{proof}

The enumeration of minimal $(k, \ell)$-tricrucial permutations now follows from the hook-length formula, \Cref{lem_for_min}, \Cref{tricrucial_rsk}, and \Cref{bijection}.

\begin{thm} \label{tricrucial_min}
Let $k \geq \ell$. The number of minimal $(k, \ell)$-tricrucial permutations is $$ (1 + \delta_{k, \ell}) \sum_{i=1}^{\ell-1} \frac{i}{2\ell-i-2} {k+i-4 \choose i-1} {2\ell-i-2 \choose \ell-1}.$$ 
\end{thm}
\begin{proof}
First, suppose that $k > \ell$. Then by \Cref{tricrucial_rsk}, there is only one filling of $P_{n}(\pi)$, and $Q_{n}(\pi)$ of the same shape as $P_{n}(\pi)$ must satisfy $q_{1, k-1} < q_{2, 2}$. Now, if $q_{1, k-1} < q_{2,2}$, then none of $q_{2, 2} < \cdots < q_{\ell-1, 2}$ can be less than $q_{1, k-1}$. So, the exact set of entries that are less than $q_{1, k-1}$ must be $\{q_{1, 1}, \ldots, q_{1, k-2}\} \cup \{q_{1,1}, \ldots, q_{i, 1}\}$ for some $i$. Now, the number of ways to fill $Q_{n}(\pi)$ is given by \Cref{lem_for_min} to be
\begin{equation} \label{tricrucial_1}
 \frac{i}{2\ell-i-2} {k+i-4 \choose i-1}{2 \ell-i-2 \choose \ell-1}.   
\end{equation}
So, we finish by summing \Cref{tricrucial_1} across $1 \leq i \leq \ell-1$. 

Next, let $k = \ell$. If $P_{n}(\pi)$ has shape $(k-1, \underbrace{2, \ldots, 2}_{k-2})$, we proceed as in the case of $k \ne \ell$ to conclude that the number of $(P_{n}(\pi), Q_{n}(\pi))$ for which $\pi$ is minimal $(k, \ell)$-tricrucial is 
\begin{equation} \label{tricrucial_2}
    \sum_{i=1}^{\ell-1} \frac{i}{2\ell-i-2} {k+i-4 \choose i-1} {2\ell-i-2 \choose \ell-1}.
\end{equation} 
Similarly, when $P_{n}(\pi)$ has shape $(k-1, k-1, \underbrace{1, \ldots, 1}_{k-3})$, the filling of $P_{n}(\pi)$ is fixed, and it follows from \Cref{bijection} that the number of ways to fill $Q_{n}(\pi)$ is also given by \Cref{tricrucial_2}. 

We thus finish by doubling \Cref{tricrucial_2}. 
\end{proof}

\subsection{Next minimal tricrucial permutations}

To count the next minimal $(k, \ell)$-tricrucial permutations, we establish the following lemma, which counts the number of $Q_{n}(\pi)$  for some next minimal $(k, \ell)$-tricrucial permutation $\pi$ that satisfy $q_{2,2} < q_{1, k-1}$. 

\begin{lem} \label{lem_for_next}
The number of standard Young tableaux $T$ of shape $(k-1, 3, \underbrace{2, \ldots, 2}_{\ell-3})$ that satisfy $t_{2, 3} < t_{3, 2}$ and contain the exact set of entries less than $t_{1, k-1}$ in $\{b_{1, 1}, b_{1,2}, \ldots, b_{1, k-2}\} \cup \{b_{1, 1}, b_{2,1}, \ldots, b_{l-i-1, 1}\} \cup \{b_{2, 2}\}$ is
$$\frac{(l+k-i-3)(k-3)(l-i-2)}{(k-2)(l+i-1)} {l+k-i-5 \choose k-3}{l+i-1 \choose l-1}.$$ 
\end{lem}
\begin{proof}

We begin by letting $S$ be the standard Young tableau formed by  $\{b_{1, 1}, b_{1,2}, \ldots, b_{1, k-1}\}$, $\{b_{1, 1}, b_{2,1}, \ldots, b_{l-i-1, 1}\},$ and $\{b_{2, 2}\}$. Let $U$ be the filled skew Young diagram of boxes formed by the boxes that are in $T$ but not in $S$. Refer to \Cref{fig:t_t'''} for the position of $S$ and $U$ in $T$. 

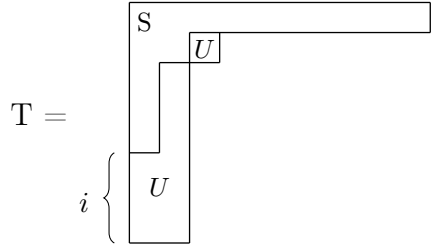
\begin{figure}[h]
\center 
\begin{tikzpicture}[scale=0.4]

\draw [line width=0.5](0,0)--(2,0);
\draw [line width=0.5](0,0)--(0,8);
\draw [line width=0.5](0,8)--(10,8);
\draw [line width=0.5](10,8)--(10,7);
\draw [line width=0.5](10,7)--(2,7);
\draw [line width=0.5](2,0)--(2,7);
\draw [line width=0.5](1,3)--(1,6);
\draw [line width=0.5](1,6)--(3,6);
\draw [line width=0.5](3,6)--(3,7);
\draw [line width=0.5](0,3)--(1,3);

\node [below] at (-3,5){{\large T = }};

 \draw [decorate,decoration={brace,amplitude=4pt, mirror},xshift=0.5cm,yshift=0pt]
      (-1,3) -- (-1,0) node [midway,right,xshift=.1cm] {};      
\node [below] at (-1.5,2){{$i$}};      
      %% brace a top
      
\node [below] at (0.5,8){{\small S}};
\node [below] at (1,2.5){{\small $U$}};
\node [below] at (2.5, 7.1){{\small $U$}};
\end{tikzpicture} 
\caption{The position of $S$ in $T$} \label{fig:t_t'''}
\end{figure}
Since the elements in $U$ are always greater than those of $S$, the partial orders between an element in $S$ and an element in $U$ given by the shape of $T$ are always respected. In addition, the partial order $t_{2, 3} < t_{3,2}$ that we must additionally satisfy from the statement of the lemma is in between two entries of $U$. It thus remains to fill $S$ and $U$ independently such that all partial orders between two elements of $S$ and between two elements of $U$ by the shape of $T$ are respected. 

We first fill $S$. We have that $t_{1, k-1} = k+\ell-i-2$, and $\{t_{1, 1}, t_{1, 2}, \ldots, t_{1, k-2}\} \cup \{t_{1,1}, \ldots, t_{\ell-i-1, 1}\} \cup \{t_{2,2}\}= \{1, 2, \ldots, k+\ell-i-3\}$. Thus, by the hook-length formula, there are \begin{equation} \label{lem_4_eq_1}
    \frac{(\ell-i-2)(k-3)}{(\ell+k-i-4)}{\ell+k-i-3 \choose k-2}
\end{equation} ways of filling $S$. 

Next, we fill $U$. Recall from the statement of the lemma that it must be that $t_{2, 3} < t_{3, 2}$. So, we first shift the boxes in $U \setminus \{b_{2,3}\}$ by one box to the right to arrive at $U'$. \Cref{fig:u_u'_relation} demonstrates how to arrive at $U'$ from $U$. Next, we take the complement of $U'$ and rotate it by 180 degrees. Then, the resulting Young diagram of shape $(\underbrace{2, \ldots, 2}_{i}, \underbrace{1, \ldots, 1}_{\ell-i-2})$ must be a standard Young tableau as its rows and columns increase.  

\begin{figure}[h]
\begin{multicols}{2}
\center 
\center 
\begin{tikzpicture}[scale=0.4]

\draw [line width=0.5](0,0)--(2,0);
\draw [line width=0.5](0,0)--(0,3);
\draw [line width=0.5](3,7)--(2,7);
\draw [line width=0.5](2,0)--(2,7);
\draw [line width=0.5](1,3)--(1,6);
\draw [line width=0.5](1,6)--(3,6);
\draw [line width=0.5](3,6)--(3,7);
\draw [line width=0.5](0,3)--(1,3);

\node [below] at (-6,5){{\large U = }};

 \draw [decorate,decoration={brace,amplitude=4pt, mirror},xshift=0.5cm,yshift=0pt]
      (-1,3) -- (-1,0) node [midway,right,xshift=.1cm] {}; 
 \draw [decorate,decoration={brace,amplitude=4pt, mirror},xshift=0.5cm,yshift=0pt]
      (-1,7) -- (-1,3) node [midway,right,xshift=.1cm] {}; 
\node [below] at (-1.5,2){{$i$}}; 
\node [below] at (-3,5.7){{$\ell-i-2$}}; 
      %% brace a top
\node[single arrow, draw=black, fill=black, 
      minimum width = 2pt, single arrow head extend=3pt,
      minimum height=10mm, scale = 0.4] at (3
    , 3) {};
\end{tikzpicture}  \columnbreak
\center 
\begin{tikzpicture}[scale=0.4]

\draw [line width=0.5](0,0)--(2,0);
\draw [line width=0.5](0,0)--(0,3);
\draw [line width=0.5](1,7)--(2,7);
\draw [line width=0.5](2,0)--(2,7);
\draw [line width=0.5](1,3)--(1,7);
\draw [line width=0.5](0,3)--(1,3);

\node [below] at (-6,5){{\large U' = }};

 \draw [decorate,decoration={brace,amplitude=4pt, mirror},xshift=0.5cm,yshift=0pt]
      (-1,3) -- (-1,0) node [midway,right,xshift=.1cm] {}; 
 \draw [decorate,decoration={brace,amplitude=4pt, mirror},xshift=0.5cm,yshift=0pt]
      (-1,7) -- (-1,3) node [midway,right,xshift=.1cm] {}; 
\node [below] at (-1.5,2){{$i$}}; 
\node [below] at (-3,5.7){{$\ell-i-2$}}; 
      %% brace a top
\end{tikzpicture}
\end{multicols} 
\caption{The construction of $U'$ from $U$} \label{fig:u_u'_relation}
\end{figure}
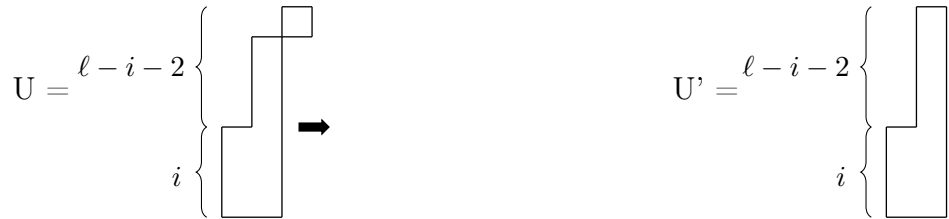

Therefore, by the hook-length formula, there are \begin{equation} \label{lem_4_eq_2}
    \frac{\ell-i-1}{\ell+i-1} {\ell+i-1 \choose i}
\end{equation} ways of filling $U$. 

We finish by multiplying \Cref{lem_4_eq_1} and \Cref{lem_4_eq_2}. \end{proof}

We next establish the following lemma, which counts the number of distinct $Q_{n}(\pi)$ for $(k, \ell)$-tricrucial permutations $\pi$ that satisfy $q_{1, k-1} < q_{2,2}$. 

\begin{lem} \label{lem_for_next''}
The number of standard Young tableaux $T$ of shape $(k-1, 3, \underbrace{2, \ldots, 2}_{\ell-3})$ that contain the exact set of entries less than $t_{1, k-1}$ in $\{b_{1, 1}, b_{1,2}, \ldots, b_{1, k-2}\} \cup \{b_{1, 1}, b_{2,1}, \ldots, b_{l-i-1, 1}\}$ and satisfy $t_{2,2} =  t_{1, k-1} + i-j +1$ is
$$(\ell-j-2){\ell+k-i-5 \choose k-3}{\ell+j-2 \choose \ell-2}.$$
\end{lem}
\begin{proof}

Let $S$ be the standard Young tableaux formed by the boxes $\{b_{1, 1}, b_{1,2}, \ldots, b_{1, k-1}\}$ and $\{b_{1, 1}, b_{2,1}, \ldots, b_{l-i-1, 1}\}$. Refer to \Cref{fig:s_t} for the position of $S$ in $T$. 

\begin{figure}[h]
\center 
\begin{tikzpicture}[scale=0.4]

\draw [line width=0.5](0,0)--(2,0);
\draw [line width=0.5](0,0)--(0,8);
\draw [line width=0.5](0,8)--(10,8);
\draw [line width=0.5](10,8)--(10,7);
\draw [line width=0.5](10,7)--(1,7);
\draw [line width=0.5](2,0)--(2,6);
\draw [line width=0.5](1,5)--(1,7);
\draw [line width=0.5](3,6)--(3,7);
\draw [line width=0.5](0,5)--(1,5);
\draw [line width=0.5](3,6)--(2,6);

\node [below] at (-3,5){{\large T = }};

 \draw [decorate,decoration={brace,amplitude=4pt, mirror},xshift=0.5cm,yshift=0pt]
      (-1,5) -- (-1,0) node [midway,right,xshift=.1cm] {};      
\node [below] at (-1.5,3){{$i$}};      
      %% brace a top
      
\node [below] at (0.5,8){{\small S}};
\node [below] at (1,3.5){{\small $T \setminus S$}};
\end{tikzpicture} 
\caption{The position of $S$ in $T$} \label{fig:s_t}
\end{figure}
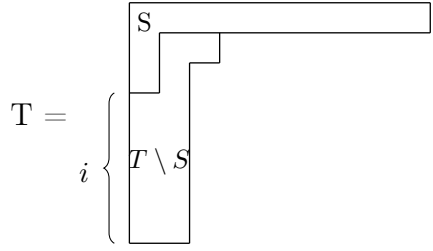

Since the elements in $T \setminus S$ are always greater than that of $S$, the partial orders given by the shape of $T$ between an element in $S$ and an element in $T \setminus S$ are always respected. It thus remains to fill $S$ and $S \setminus T$ independently such that all partial orders between two elements of $S$ and between two elements of $T \setminus S$ given by the shape of $T$ are respected. 

We first fill $S$. We have that $t_{1, k-1} = k+\ell-i-3$, and $\{t_{1, 1}, t_{1, 2}, \ldots, t_{1, k-2}\} \cup \{t_{1,1}, \ldots, t_{\ell-i-1, 1}\}= \{1, 2, \ldots, k+\ell-i-4\}$. Thus, by the hook-length formula, there are \begin{equation} \label{lem_5_eq_1}
    {\ell+k-i-5 \choose k-3}
\end{equation} ways of filling $S$.

Next, we fill $T \setminus S$. Recall that $t_{2,2} =t_{1, k-1}+i-j+1$. First, because the only entries of $T\setminus S$ that can be less than $t_{2,2}$ are in row $1$, it must be that
\begin{align*}
    t_{\ell-i} = t_{1, k-1} & +1 = k+\ell-i-2, \\ t_{\ell-i+1} = t_{1, k-1}& +2 = k+\ell-i-1, \\  &\vdots \\ t_{\ell-j-1} = t_{1,k-1} &+ i-j = k+\ell-j-1.
\end{align*}

 Next, because $t_{2,2} = k+\ell-j-2 < t_{2,3}$, there are 
\begin{equation} \label{lem_5_eq_3}
    (k + 2\ell-4)-(k+\ell-j-2) = \ell+j-2
\end{equation} ways of filling $t_{2,3}$. 

Let $U$ be the filled skew Young diagram formed by the boxes $\{b_{\ell-j, 1}, b_{\ell-i+j+1, 1}, \ldots, b_{\ell-1, 1}\} \cup \{b_{3,2}, b_{4,2}, \ldots, b_{\ell-1,2}\}$. Refer to \Cref{fig:u_t} for the position of $U$ in $T$. It remains to fill $U$ with $\{k+\ell-j-1,k+\ell-j, \ldots, k+2 \ell -4\} \setminus \{t_{3,2}\}$. Now, add $1$ to every entry of $U$ that is less than $t_{3,2}$, take the complement, and rotate by $180$ degrees. Then we must have a standard Young tableau of shape $(\underbrace{2, \ldots, 2}_{j}, \underbrace{1, \ldots, 1}_{\ell-j-2})$ as each row and column must be increasing. 

\begin{figure}[h]
\center 
\begin{tikzpicture}[scale=0.4]

\draw [line width=0.5](0,0)--(2,0);
\draw [line width=0.5](0,0)--(0,8);
\draw [line width=0.5](0,8)--(10,8);
\draw [line width=0.5](10,8)--(10,7);
\draw [line width=0.5](10,7)--(3,7);
\draw [line width=0.5](2,0)--(2,6);
\draw [line width=0.5](1,3)--(1,6);
\draw [line width=0.5](3,6)--(3,7);
\draw [line width=0.5](0,3)--(1,3);
\draw [line width=0.5](3,6)--(1,6);

\node [below] at (-3,5){{\large T = }};

 \draw [decorate,decoration={brace,amplitude=4pt},xshift=0.5cm,yshift=0pt]
      (2,6) -- (2,0) node [midway,right,xshift=.1cm] {}; 

 \draw [decorate,decoration={brace,amplitude=4pt, mirror},xshift=0.5cm,yshift=0pt]
      (-1,3) -- (-1,0) node [midway,right,xshift=.1cm] {};      
\node [below] at (-2.5,2.2){{$j$}};   \node [below] at (4.2,3.5){{$\ell-3$}};
      %% brace a top
      
\node [below] at (1,2.1){{\small $U$}};
\end{tikzpicture} 
\caption{The position of $U$ in $T$} \label{fig:u_t}
\end{figure}
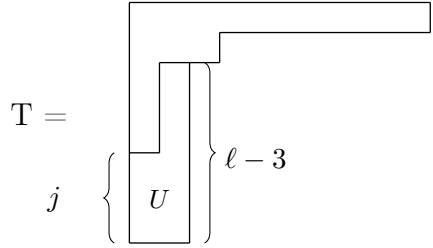

Therefore, by the hook-length formula, there are \begin{equation} \label{lem_5_eq_2}
    \frac{\ell-j-2}{\ell+j-2} {\ell+j-2 \choose \ell-2}
\end{equation} ways of filling $U$. 

We finish by multiplying \Cref{lem_5_eq_1}, \Cref{lem_5_eq_3}, and \Cref{lem_5_eq_2}.
\end{proof}

As we proceed, we let $\mathrm{rev}(\pi)$ to be the permutation $\pi_{n} \pi_{n-1} \cdots \pi_{1}$, and the $(k, \ell)$-left(respectively, top)-crucial permutations to be the $(k, \ell)$-crucial permutations for which the only allowed extensions are to the left (respectively, from above). The following lemma helps us count the $(k, k)$-tricrucial permutations. 

\begin{lem} \label{next_bijection}
The permutation $\pi$ is a $(k,k)$-tricrucial permutation if and only if $\mathrm{rev}(\pi)$ is a $(k, k)$-tricrucial permutation. 
\end{lem}
\begin{proof}
An increasing (respectively, decreasing) subsequence of length $k$ in $\pi$ is a decreasing (respectively, increasing) subsequence of length $k$ in $\mathrm{rev}(\pi)$. Therefore, $\pi$ is left (respectively, right, top)-crucial if and only if $\mathrm{rev}(\pi)$ is right (respectively, left, top)-crucial as sought. 
\end{proof}

We proceed to count the next minimal $(k, \ell)$-tricrucial permutations. Our enumeration follows from the hook-length formula, \Cref{lem_for_next_3}, \Cref{lem_next_transpose}, \Cref{lem_next_one}, \Cref{lem_for_min}, \Cref{lem_for_next}, \Cref{lem_for_next''}, and \Cref{next_bijection}. As in the proof of \Cref{bicrucial_next}, we deal with the cases $(k, \ell) \in \{(4, 3), (4,4)\}$ separately, because unlike $(k, \ell)$ with larger $k$ for which $k = \ell$ or $k = \ell+1$, the shape of $P_{n}(\pi)$ (and $Q_{n}(\pi)$) are fixed when $(k, \ell) \in \{(4, 3), (4,4)\}$. 

\begin{thm} \label{next_min_tri}
Suppose that $k \geq \ell$. The number of next minimal $(4,3)$-tricrucial permutations is $25$, the number of next minimal $(4,4)$-tricrucial permutations is $540$, and the number of next minimal $(k, \ell)$-tricrucial permutations for $k \geq \ell$ and $(k, \ell) \not\in \{(4,3), (4,4)\}$ is
\begin{tiny}
\begin{align*}
     (1+ \delta_{k, \ell})(k+2\ell-7) \sum_{i=0}^{\ell-2} {\ell+k-i-5 \choose k-3} \left( \sum_{j=0}^{i+1} (\ell-j-2) {\ell+j-2 \choose \ell-2} + \frac{(\ell+k-i-3)(k-3)(\ell-i-2)}{(k-2)(\ell+i-1)} {\ell+i-1 \choose \ell-1} \right)  \\ +  \delta_{k-1, \ell} \sum_{i=1}^{k-1} \frac{i}{2k-i-2} {\ell+i-4 \choose i-1} {2k-i-2 \choose k-1},
\end{align*}
\end{tiny}
\end{thm}
\begin{proof}
By Theorem 31 of Avgustinovich, Kitaev, and Taranenko \cite{akt}, the length of a minimal $(k, \ell)$-tricrucial permutation is $k+2 \ell-5$. Therefore, the length of a next minimal $(k, \ell)$-tricrucial permutations is $k + 2 \ell-4$. In addition, $P_{n}(\pi)$ must have $k-1$ columns and $\ell-1$ rows by conditions $1$ and $2$ of \Cref{tricrucial_1}.

First, suppose that $p_{\ell-1, 1} < p_{1, k-1}$. Then, to satisfy condition 3 of \Cref{tricrucial_char}, column $2$ of $P_{n}(\pi)$ must have $\ell-1$ boxes. Thus, for $P_{n}(\pi)$ to have $k-1$ boxes in row $1$, $\ell-1$ boxes in columns $1$ and $2$, and $k+2\ell-4$ boxes total, $P_{n}(\pi)$ must have the shape $(k-1, 3, \underbrace{2, \ldots, 2}_{\ell-3})$. Now, by \Cref{lem_for_next_3}, there are 
\begin{equation} \label{next_min_tri_eq_1}
    k + 2 \ell - 7
\end{equation} ways to fill $P_{n}(\pi)$ such that it satisfies conditions $3$ and $4$ of \Cref{tricrucial_char}. 

Now, from \Cref{tricrucial_char}, $Q_{n}(\pi)$ must be a standard Young tableau with the same shape as $P_{n}(\pi)$ that contains an increasing sequence of elements $q_{1, k-1}, q_{2, m_2}, \ldots, q_{\ell-1, m_{\ell-1}}$ for some $m_2, \ldots, m_{\ell-1}$. We divide cases based on whether $q_{2,2} < q_{1, k-1}$ or not. 

Now, suppose that $q_{2,2} < q_{1, k-1}$. Then, for $Q_{n}(\pi)$ to contain an increasing sequence of elements $q_{1, k-1}, q_{2, m_2}, \ldots, q_{\ell-1, m_{\ell-1}}$ for some $m_2, \ldots, m_{\ell-1}$, it must be that $q_{1, k-1}< q_{2,3}< q_{3,2}$. Therefore, by \Cref{lem_for_next}, the number of $Q_{n}(\pi)$ in this case is given by summing 
\begin{equation} \label{next_min_tri_eq_2}
    \frac{(l+k-i-3)(k-3)(l-i-2)}{(k-2)(l+i-1)} {l+k-i-5 \choose k-3}{l+i-1 \choose l-1}
\end{equation} across $0 \leq i \leq \ell-2$. 

By the RSK correspondence, we finish this case by summing \Cref{next_min_tri_eq_2} across $0 \leq i \leq \ell-2$ and multiplying it with \Cref{next_min_tri_eq_1}. 

Now, suppose that $q_{1, k-1} < q_{2,2}$. Then, by \Cref{lem_for_next''}, the number of ways to fill $Q_{n}(\pi)$ is given by summing
\begin{equation} \label{next_min_tri_eq_3}
    (\ell-j-2){\ell+k-i-5 \choose k-3}{\ell+j-2 \choose \ell-2}
\end{equation} across $0 \leq i \leq \ell-2$ and $0 \leq j \leq i+1$. 

By the RSK correspondence, we finish this case by summing \Cref{next_min_tri_eq_2} across $0 \leq i \leq \ell-2$ and $0 \leq j \leq i+1$ and multiplying it with \Cref{next_min_tri_eq_3}. 

Now, suppose that $p_{1, k-1} < p_{\ell-1, 1}$. Then, because all entries in row $1$ are less than $p_{\ell-1, 1}$, for condition $4$ of \Cref{tricrucial_char} to hold, row $2$ must have all $k-1$ boxes. Therefore, the minimum number of boxes in $P_{n}(\pi)$ is $2k + \ell-5$. Since $\pi$ has $k + 2 \ell- 4$, it must be that $k \in \{\ell, \ell+1\}$. 

First, suppose that $k = \ell$. Then, for $P_{n}(\pi)$ to have $k-1$ boxes in rows $1$ and $2$, $k-1$ boxes in column $1$, and $3k-4$ boxes total, $P_{n}(\pi)$ must have the shape $(k-1, k-1, 2, \underbrace{1, \ldots, 1}_{k-4})$. Now, as proved in \cite{enum}, the shape of the RSK correspondence of $\pi$ must be the transpose of that of $\mathrm{rev}(\pi)$. Furthermore, by \Cref{next_bijection}, $\pi$ is a $(k, k)$-tricrucial permutation if and only if $\mathrm{rev}(\pi)$ is. Therefore, the number of $(k, \ell)$-tricrucial permutations must be the same as in the case of $p_{\ell-1, 1} < p_{1, k-1}$. Namely, it is given by: 

\begin{align}
    & (k+2\ell-7) \sum_{i=0}^{\ell-2} {\ell+k-i-5 \choose k-3} \\ & \left( \sum_{j=0}^{i+1} (\ell-j-2) {\ell+j-2 \choose \ell-2} + \frac{(\ell+k-i-3)(k-3)(\ell-i-2)}{(k-2)(\ell+i-1)} {\ell+i-1 \choose \ell-1} \right) \nonumber.
\end{align}

Now, if $k = \ell+1$, then for $P_{n}(\pi)$ to have $k-1$ boxes in rows $1$ and $2$, $\ell-1$ boxes in column $1$, and $k+2\ell-4$ boxes total, $P_{n}(\pi)$ must have shape $(k-1, k-1, \underbrace{1, \ldots, 1}_{\ell-3})$. Therefore, by \Cref{lem_next_one}, there is exactly one way to fill $P_{n}(\pi)$. 

Now, we fill $Q_{n}(\pi)$. By \Cref{tricrucial_char}, $Q_{n}(\pi)$ must have the same shape as $P_{n}(\pi)$ and contain an increasing sequence of elements $q_{1, k-1}, q_{2, m_2}, \ldots, q_{\ell-1, m_{\ell-1}}$ for some $m_2, \ldots, m_{\ell-1}$. By \Cref{bijection}, the number of ways to fill $Q_{n}(\pi)$ is the same as the number of standard Young tableau $Q'$ of shape $(\ell-1, \underbrace{2, \ldots, 2}_{k-2})$ that satisfies $q'_{1, k-1} < q'_{2,2}$. Therefore, by \Cref{lem_for_min}, the number of ways to fill $Q_{n}(\pi)$ is given by summing 
\begin{equation}\label{next_min_tri_eq_4}
    \frac{i}{2\ell-i-2} {k+i-4 \choose i-1}{2 \ell-i-2 \choose \ell-1}
\end{equation} across $1 \leq i \leq k-1$. 

By the RSK correspondence, we finish this case by summing \Cref{next_min_tri_eq_4} across $1 \leq i \leq k-1$.  

To complete the enumeration of the next minimal $(k, \ell)$-tricrucial permutations, we verify via computer search that the number of next minimal $(4,3)$-tricrucial permutations is $25$ and the number of next minimal $(4,4)$-tricrucial permutations is $540$.
\end{proof}

\subsection{$(k, 3)$-tricrucial permutations}

To count the $(k, 3)$-tricrucial permutations, we first establish the following lemma, which counts the number of $P_{n}(\pi)$ for some $(k, 3)$-tricrucial permutation $\pi$ when $k +1 \leq n \leq 2k-3$. 

\begin{lem} \label{lem_for_3}
Let $k+1 \leq n \leq 2k-3$. The number of standard Young tableaux $T$ of shape $(k-1, n-k+1)$ that satisfy $t_{2, n-k+1}=n$ and contain an increasing sequence of elements $t_{2, 1}, t_{i_2, 2}, \ldots, t_{i_{k-1}, k-1}$ for some $i_2, i_3, \ldots, i_{k-1}$ is $$\frac{2k-n}{n} {n \choose k} - \frac{2k-n-1}{n-1}{n-1 \choose k-1}.$$
\end{lem}
\begin{proof}
With $t_{2, n-k+1} = n$ fixed, the number of standard Young tableaux on the remaining shape $(k-1, n-k)$ is given by the hook-length formula to be \begin{equation} \label{lem_1_eq_1}
    \frac{2k-n}{n} {n \choose k}.
\end{equation}
Next, we claim that if a standard Young tableau of shape $(k-1, n-k+1)$ satisfies $t_{2, n-k+1} = n$, an increasing sequence of elements $t_{2, 1}, t_{i_2, 2}, \ldots, t_{i_{k-1}, k-1}$ for some $i_2, i_3, \ldots, i_{k-1}$ does not exist if and only if $t_{2, j} > t_{1, j+1}$ for all $1 \leq j \leq n-k+1$. We first show the condition $t_{2, j} > t_{1, j+1}$ for all $1 \leq j \leq n-k+1$ is necessary. Assume that such a sequence exists. Then the condition $t_{2, j} > t_{1, j+1}$ for all $1 \leq j \leq n-k+1$ forces $i_2 = \cdots = i_{n-k+1} = 2$. However, since $t_{2, n-k+1} = n$, we cannot have $n=t_{2, n-k+1} < t_{1, n-k+2}$. Therefore, such a sequence does not exist. We now show that the condition $t_{2, j} > t_{1, j+1}$ for all $1 \leq j \leq n-k+1$ is sufficient. If $t_{2, j} < t_{1, j+1}$ for some $j$ such that $1 \leq j \leq n-k+1$, then let $i_2 = \cdots = i_{j} = 2$ and $i_{j+1} = \cdots = i_{k-1} = 1$. Then the sequence $t_{2, 1}, t_{2, 2}, \ldots, t_{2, j}, t_{1, j+1}, t_{1, j+2} \ldots, t_{1, k-1}$ increases, since each row of a standard Young tableau increases and $t_{2, j} < t_{1, j+1}$. 

Now, we enumerate the standard Young tableaux $T$ of shape $(k-1, n-k+1)$ that satisfy $t_{2, n-k+1} = n$ and $t_{2, j} > t_{1, j+1}$ for all $1 \leq j \leq n-k+1$. To do so, we first shift the boxes in row $2$ by one box to the right. We then eliminate the box that contains $t_{1,1}=1$ and subtract $1$ from each entry of $T$. Then the resulting Young diagram of shape $(k-2, n-k+1)$ is a standard Young tableau $T'$ that satisfies $t'_{2, n-k+1} = n-1$: the rows in $T'$ increase as the rows in $T$ increase, and the columns in $T'$ increase as $t_{2, j} > t_{1, j+1}$ for all $1 \leq j \leq n-k+1$. \Cref{t_t'_relation} demonstrates how to arrive at $T'$ from $T$. 

\begin{figure}[h]
\begin{multicols}{2} \center 
\begin{tikzpicture}[scale=0.35]

\draw [line width=0.5](0,1)--(11,1);
\draw [line width=0.5](0,2)--(11,2);
\draw [line width=0.5](0,0)--(0,2);
\draw [line width=0.5](1,0)--(1,2);
\draw [line width=0.5](2,0)--(2,2);

\draw [line width=0.5](5,0)--(5,2);
\draw [line width=0.5](6,0)--(6,2);

\draw [line width=0.5](7,1)--(7,2);
\draw [line width=0.5](0,0)--(6,0);

\draw [line width=0.5](10,1)--(10,2);
\draw [line width=0.5](11,1)--(11,2);

\draw [decorate,decoration={brace,amplitude=5pt,mirror,raise=4ex}]
  (0,1) -- (6,1) node[midway,yshift=-3em]{$n-k+1$};
  
\draw [decorate,decoration={brace,amplitude=5pt, raise=4ex}]
  (0,1) -- (11,1) node[midway,yshift=-3em]{};

\node [below] at (-1.3,1.5){{\large $T = $}};
  
\node [below] at (5.5,4.7){{$k-1$}};

\node [below] at (3.5,1.9){{\huge $\ldots$}};
\node [below] at (5.3,1.1){{ $n$}};
\node [below] at (0.35,2.2){{ $1$}};
\node [below] at (8.5,1.9){{\huge $\ldots$}};
\node [below] at (3.5,0.9){{\huge $\ldots$}};
\node[single arrow, draw=black, fill=black, 
      minimum width = 2pt, single arrow head extend=3pt,
      minimum height=10mm, scale = 0.4] at (7
    , 0.5) {};% length of arrow

\end{tikzpicture}
\columnbreak \center 
\begin{tikzpicture}[scale=0.35]

\node [below] at (-1.3,1.5){{\large $T' = $}};

\draw [line width=0.5](0,1)--(12,1);
\draw [line width=0.5](0,2)--(12,2);
\draw [line width=0.5](0,1)--(0,2);
\draw [line width=0.5](1,0)--(1,2);
\draw [line width=0.5](2,0)--(2,2);

\draw [line width=0.5](5,0)--(5,2);
\draw [line width=0.5](6,0)--(6,2);
\draw [line width=0.5](8,0)--(8,2);

\draw [line width=0.5](1,0)--(8,0);

\draw [line width=0.5](11,1)--(11,2);
\draw [line width=0.5](12,1)--(12,2);

\node [below] at (3.5,1.9){{\huge $\ldots$}};
\node [below] at (0.5,2.2){{\Cross}};
\node [below] at (7,1.1){{\tiny $n-1$}};
\node [below] at (0.35,2.2){{ $1$}};
\node [below] at (9.5,1.9){{\huge $\ldots$}};
\node [below] at (3.5,0.9){{\huge $\ldots$}};
\draw [decorate,decoration={brace,amplitude=5pt,mirror,raise=4ex}]
  (1,1) -- (8,1) node[midway,yshift=-3em]{$n-k+1$};
  
\draw [decorate,decoration={brace,amplitude=5pt, raise=4ex}]
  (1,1) -- (12,1) node[midway,yshift=-3em]{};
  
\node [below] at (6.3,4.7){{$k-2$}};

\end{tikzpicture}
\end{multicols}
    \caption{The construction of $T'$ from $T$}
    \label{t_t'_relation}
\end{figure}
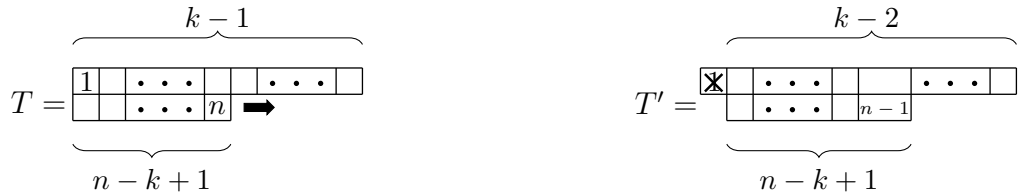

Thus, our count reduces to counting the number of standard Young tableaux $T'$ of shape $(k-2, n-k+1)$ that satisfy $t'_{2, n-k} = n-1$. After removing the box $t'_{2, n-k} = n-1$, our count further reduces to counting the number of standard Young tableaux with the shape $(k-2, n-k)$. Therefore, the number of standard Young tableaux $T$ of shape $(k-1, n-k+1)$ that satisfy $t_{2, n-k+1}=n$ and do not contain an increasing sequence of elements $t_{2, 1}, t_{i_2, 2}, \ldots, t_{i_{k-1}, k-1}$ for some $i_2, i_3, \ldots, i_{k-1}$ is given by the hook-length formula to be \begin{equation} \label{lem_1_eq_2}
    \frac{2k-n-1}{n-1} {n-1 \choose k-1}. 
\end{equation}
We finish by subtracting \Cref{lem_1_eq_2} from \Cref{lem_1_eq_1}.
\end{proof}

We now proceed to our enumeration of the $(k, 3)$-tricrucial permutations for $k \geq 3$. Our enumeration is a culmination of the hook-length formula, \Cref{tricrucial_char}, and \Cref{lem_for_3}. 

\begin{thm} \label{tricrucial_3}
Let $k \geq 3$. The number of $(k, 3)$-tricrucial permutations of length $n$ for $k+1 \leq n \leq 2k-3$ is $$\frac{2k-n}{n} {n \choose k} \left( \frac{2k-n}{n}{n \choose k} - \frac{2k-n-1}{n-1} {n-1 \choose k-1} \right).$$ For $n = 2(k-1)$, it is $$ \left( \frac{1}{k} {2k-2 \choose k-1} \right)^2.$$  
\end{thm}
\begin{proof}
 Since $\ell = 3$, the third condition of \Cref{tricrucial_char} is equivalent to $p_{1, k-1} < p_{2, n-k+1}$. Because the rows of $P_{n}(\pi)$ increase, either $p_{1, k-1} = n$ or $p_{2, n-k+1} = n$. Therefore, the third condition of \Cref{tricrucial_char} is equivalent to $p_{2, n-k+1}=n$. Similarly, the fifth condition of \Cref{tricrucial_char} is equivalent to $q_{2, n-k+1} = n$. 

Let $k+1  \leq n \leq 2k-3$. Then from \Cref{lem_for_3}, the number of possible fillings of $P_{n}(\pi)$ is given by  
\begin{equation} \label{tricrucial_4}
      \frac{2k-n}{n}{n \choose k} - \frac{2k-n-1}{n-1} {n-1 \choose k-1}.
\end{equation} 
Now, we fill $Q_{n}(\pi)$. Since $q_{2, n-k+1}=n$, we must fill the remaining tableau of shape $(k-1, n-k)$. Thus, from the hook-length formula, the number of possible fillings of $Q_{n}(\pi)$ is given by
\begin{equation} \label{tricrucial_5}
    \frac{2k-n}{n} {n \choose k}.
\end{equation} 

By the RSK correspondence, we finish by multiplying \Cref{tricrucial_4} and \Cref{tricrucial_5}.

Now, suppose that $n = 2(k-1)$. Then the third condition of \Cref{tricrucial_char} is always met, because $p_{1, k-1} < p_{2, k-1} = n$. Likewise, the fifth condition is always met. Similarly, the fourth condition of \Cref{tricrucial_char} is also always met, because $p_{2, 1} < p_{2, 2} < \cdots < p_{2, k-1}$. Thus, the only restrictions on $P_{n}(\pi)$ and $Q_{n}(\pi)$ are that they are standard Young tableaux of shape $(k-1, k-1)$. Therefore, by the hook-length formula, there are \begin{equation} \label{tricrucial_eq_3}
    \frac{1}{k}{2k-2 \choose k-1}
\end{equation} possible fillings for each of $P_{n}(\pi)$ and $Q_{n}(\pi)$. 

By the RSK correspondence, we finish by squaring \Cref{tricrucial_eq_3}. 
\end{proof}

That all $(k, \ell)$-tricrucial permutations for $k \geq \ell$ must have length at least $k+2\ell-5$ and at most $(k-1)(\ell-1)-1$ follows from Theorem 32 of Avgustinovich, Kitaev, and Taranenko \cite{akt}. Therefore, \Cref{tricrucial_3} gives a complete enumeration of $(k, 3)$-tricrucial permutations. 

\section{Quadrocrucial Permutations}
\label{quadro}
Our next result is that we enumerate the next minimal $(k, \ell)$-quadrocrucial permutations. We begin by citing the characterization of $(P_{n}(\pi), Q_{n}(\pi))$ when $\pi$ is quadrocrucial as presented by Avgustinovich, Kitaev, and Taranenko \cite{akt}. 

\begin{thm}[Avgustinovich, Kitaev, and Taranenko \cite{akt}] \label{quadrocrucial_char}
The permutation $\pi \in S_{n}$ is $(k,\ell)$-quadrocrucial if and only if:
\begin{enumerate}
    \item the number of columns in $P_{n}(\pi)$ (and $Q_{n}(\pi)$) is $k-1$;
    \item the number of rows in $P_{n}(\pi)$ (and $Q_{n}(\pi)$) is $\ell-1$;
    \item $P_{n}(\pi)$ contains an increasing sequence of elements $p_{1, k-1}, p_{2, j_2}, \ldots, p_{\ell-1, j_{\ell-1}}$ for some $j_2, \ldots,$ $j_{\ell-1}$; 
    \item $P_{n}(\pi)$ contains an increasing sequence of elements $p_{\ell-1,1}, p_{i_2, 2}, \ldots,p_{i_{k-1}, k-1}$ for some $i_2, \ldots,$ $i_{k-1}$;
    \item $Q_{n}(\pi)$ contains an increasing sequence of elements $q_{1, k-1}, q_{2, m_2}, \ldots, q_{\ell-1, m_{\ell-1}}$ for some $m_2,$ $\ldots,$ $m_{\ell-1}$;
    \item $Q_{n}(\pi)$ contains an increasing sequence of elements $q_{\ell-1,1}, q_{r_2, 2}, \ldots,q_{r_{k-1}, k-1}$ for some $r_2, \ldots,$ $r_{k-1}$.
\end{enumerate}
\end{thm}

\subsection{Next minimal quadrocrucial permutations}

The enumeration of next minimal $(k, \ell)$-quadrocrucial permutations follows from \Cref{lem_for_next_3} and \Cref{lem_next_one}. As in the proof of \Cref{bicrucial_next}, we deal with the cases $(k, \ell) \in \{(4, 3), (4,4)\}$ separately, because unlike $(k, \ell)$ with larger $k$ for which $k = \ell$ or $k = \ell+1$, the shape of $P_{n}(\pi)$ (and $Q_{n}(\pi)$) are fixed when $(k, \ell) \in \{(4, 3), (4,4)\}$. 

\begin{thm} \label{next_min_quadro}
Suppose that $k \geq \ell$. The number of next minimal $(4,3)$-quadrocrucial permutations is $25$, the number of next minimal $(4,4)$-quadrocrucial permutations is $756$, and the number of next minimal $(k, \ell)$-quadrocrucial permutations for $k \geq \ell$ and $(k, \ell) \not\in \{(4,3), (4,4)\}$ is
\begin{equation*}
    (1+\delta_{k, \ell})(k+2\ell-7)^{2} + \delta_{k-1, \ell}.
\end{equation*} 
\end{thm}
\begin{proof}
By Theorem 36 of Avgustinovich, Kitaev, and Taranenko \cite{akt}, the length of a minimal $(k, \ell)$-quadrocrucial permutation is $k+2 \ell-5$. Therefore, the length of a next minimal $(k, \ell)$-quadrocrucial permutations is $k + 2 \ell-4$. By conditions $1$ and $2$ of \Cref{quadrocrucial_char}, $P_{n}(\pi)$ must have $k-1$ columns and $\ell-1$ rows. 

First, suppose that $p_{\ell-1, 1} < p_{1, k-1}$. Then, to satisfy condition 3 of \Cref{quadrocrucial_char}, column $2$ of $P_{n}(\pi)$ must have $\ell-1$ boxes. Thus, for $P_{n}(\pi)$ to have $k-1$ boxes in row $1$, $\ell-1$ boxes in columns $1$ and $2$, and $k+2\ell-4$ boxes total, $P_{n}(\pi)$ must have the shape $(k-1, 3, \underbrace{2, \ldots, 2}_{\ell-3})$. Now, by \Cref{lem_for_next_3}, there are 
\begin{equation} \label{next_min_quadro_eq_4}
    k + 2 \ell - 7
\end{equation} ways to fill $P_{n}(\pi)$ such that it satisfies conditions $3$ and $4$ of \Cref{quadrocrucial_char}. The same applies to $Q_{n}(\pi)$.

By the RSK correspondence, we finish this case by squaring \Cref{next_min_quadro_eq_4}.

Now, suppose that $p_{1, k-1} < p_{\ell-1, 1}$. Then because all entries in row $1$ are less than $p_{\ell-1, 1}$, for condition $4$ of \Cref{quadrocrucial_char} to hold, row $2$ must have all $k-1$ boxes. Therefore, the minimum number of boxes in $P_{n}(\pi)$ is $2k + \ell-5$. Since $\pi$ has $k + 2 \ell- 4$, it must be that $k \in \{\ell, \ell+1\}$. 

Suppose that $k = \ell$. Then, for $P_{n}(\pi)$ to have $k-1$ boxes in rows $1$ and $2$, $k-1$ boxes in column $1$, and $3k-4$ boxes total, $P_{n}(\pi)$ must have the shape $(k-1, k-1, 2, \underbrace{1, \ldots, 1}_{k-4})$. Thus, by \Cref{lem_next_transpose}, there must be 
\begin{equation} \label{next_min_quadro_eq_1}
    3k-7    
\end{equation} ways of filling $P_{n}(\pi)$. The same applies to $Q_{n}(\pi)$. 

By the RSK correspondence, we finish this case by squaring \Cref{next_min_quadro_eq_1}. 

Now, if $k = \ell+1$, for $P_{n}(\pi)$ to have $k-1$ boxes in rows $1$ and $2$, $\ell-1$ boxes in column $1$, and $k+2\ell-4$ boxes total, $P_{n}(\pi)$ must have shape $(k-1, k-1, \underbrace{1, \ldots, 1}_{\ell-3})$. Therefore, by \Cref{lem_next_one}, there is exactly one way to fill $P_{n}(\pi)$. The same applies to $Q_{n}(\pi)$. 

By the RSK correspondence, there is one additional $(k, \ell)$-quadrocrucial permutation in this case.

To complete the enumeration of the next minimal $(k, \ell)$-quadrocrucial permutations, we verify via computer search that the number of next minimal $(4,3)$-quadrocrucial permutations is $25$ and the number of next minimal $(4,4)$-quadrocrucial permutations is $756$.
\end{proof}

\subsection{$(k,3)$-quadrocrucial permutations}

Our final result is that we count the $(k, 3)$-quadrocrucial permutations for $k \geq 3$. The enumeration follows from the hook-length formula, \Cref{lem_for_min}, and \Cref{quadrocrucial_3}. 

\begin{thm} \label{quadrocrucial_3}
Let $k \geq 3$. The number of $(k, 3)$-quadrocrucial permutations of length $n$ for $k+1 \leq n \leq 2k-3$ is $$\left( \frac{2k-n}{n}{n \choose k} - \frac{2k-n-1}{n-1} {n-1 \choose k-1} \right)^2.$$ For $n = 2(k-1)$, it is $$\left( \frac{1}{k} {2k-2 \choose k-1} \right)^2.$$ 
\end{thm} 
\begin{proof}
Since $\ell = 3$, the third condition of \Cref{quadrocrucial_char} is equivalent to $p_{1, k-1} < p_{2, n-k+1}$. Since the rows of $P_{n}(\pi)$ increase, either $p_{1, k-1} = n$ or $p_{2, n-k+1} = n$. Therefore, the third condition of \Cref{quadrocrucial_char} is equivalent to $p_{2, n-k+1}=n$. Similarly, the fifth condition of \Cref{quadrocrucial_char} is equivalent to $q_{2, n-k+1} = n$. 

Now, let $k+1 \leq n \leq 2k-3$. The number of possible fillings for each of $P_{n}(\pi)$ and $Q_{n}(\pi)$ is given by \Cref{lem_for_3} to be
\begin{equation} \label{quadrocrucial_1}
     \frac{2k-n}{n}{n \choose k} - \frac{2k-n-1}{n-1} {n-1 \choose k-1}.
\end{equation} 
By the RSK correspondence, we finish by squaring \Cref{quadrocrucial_1}. 

Now, suppose that $n = 2(k-1)$. Then the third condition of \Cref{tricrucial_char} is always met, because $p_{1, k-1} < p_{2, k-1} = n$. Likewise, the fifth condition is always met. Similarly, the fourth condition of \Cref{tricrucial_char} is also always met, because $p_{2, 1} < p_{2, 2} < \cdots< p_{2, k-1}$. Likewise, the sixth condition is always met. Thus, the only restrictions on $P_{n}(\pi)$ and $Q_{n}(\pi)$ are that they are standard Young tableaux of shape $(k-1, k-1)$. Therefore, by the hook-length formula, there are \begin{equation} \label{quadrocrucial_eq_3}
    \frac{1}{k}{2k-2 \choose k-1}
\end{equation} possible fillings for each of $P_{n}(\pi)$ and $Q_{n}(\pi)$. 
By the RSK correspondence, we finish by squaring \Cref{quadrocrucial_eq_3}. 
\end{proof}

That all $(k, \ell)$-quadrocrucial permutations for $k \geq \ell$ must have length at least $k+2\ell-5$ and at most $(k-1)(\ell-1)-1$ follows from Theorem 36 of Avgustinovich, Kitaev, and Taranenko \cite{akt}. Therefore, \Cref{quadrocrucial_3} gives a complete enumeration of $(k, 3)$-quadrocrucial permutations.

\section*{Acknowledgements}
This research was conducted as a part of the 2022 University of Minnesota Duluth REU, supported by Jane Street Capital, the NSA (grant number H98230-22-1-0015), and the Harvard Herchel Smith Undergraduate Science Research Program. The author thanks Emily Zhang, Amanda Burcroff, and Joe Gallian for their thoughtful comments on the presentation of the paper.

\medskip
\begin{footnotesize}
  Y. Choi, \textsc{Harvard University,
    Cambridge, MA 02138}\par\nopagebreak
  \textit{E-mail address:} \texttt{ychoi@college.harvard.edu}
\end{footnotesize}
  
\end{document}